\documentclass[twocolumn]{revtex4-2}

\usepackage{amscd}
\usepackage{amsmath}
\usepackage{amssymb}
\usepackage{graphicx}%

\usepackage{scalerel}
\usepackage{dsfont}
\usepackage{bm}
\usepackage{float}
\usepackage{placeins}

\newcommand{\overbar}[1]{\mkern 1.5mu\overline{\mkern-3.2mu#1\mkern-1mu}\mkern 1.5mu}

\def\bi{\begin{itemize}}
\def\ei{\end{itemize}}

\usepackage{newtxtext}
\usepackage{booktabs,caption} 
\usepackage{wrapfig,overpic,setspace}

\usepackage[colorlinks=true, pdfstartview=FitV, linkcolor=blue,
            citecolor=blue, urlcolor=blue]{hyperref} 
\usepackage{subcaption}
\usepackage{sidecap}
\usepackage{xcolor}
\newcommand\wh[1]{\hstretch{2}{\hat{\hstretch{.6}{#1}}}}

\newcommand{\norm}[1]{\left\lVert#1\right\rVert}

\newtheorem{thm}{Theorem}[section]
\newtheorem{rem}{Remark}[section]
\def\bt{\begin{thm}}
\def\et{\end{thm}}
\def\br{\begin{rem}}
\def\er{\end{rem}}
\newtheorem{lem}{Lemma}[section]

\def\bea{\begin{equation} \begin{aligned}}
\def\eea{\end{aligned} \end{equation}}
\def\beas{\begin{equation*} \begin{aligned}}
\def\eeas{\end{aligned} \end{equation*}}
\def\bes{\begin{equation*}}
\def\ees{\end{equation*}}

\def\d{\, \mathrm{d}}
\def\be{\begin{equation}}
\def\ee{\end{equation}}
\def\adots{
  \mathinner{\mkern1mu\raise1pt\hbox{.}\mkern2mu\raise4pt\hbox{.}
  \mkern2mu\raise7pt\vbox{\kern7pt\hbox{.}}\mkern1mu}}
\def\bftau{{\boldsymbol{\tau}}}
\def\bflambda{{\boldsymbol{\lambda}}}
\def\s{\mathfrak{s}}
\def \c{\mathfrak{c}}

\graphicspath{{./Figures/}}

\begin{document}

\title{Optimal Parameterizing Manifolds for Anticipating Tipping Points and Higher-order Critical Transitions}

\author{Micka{\"e}l D. Chekroun}
%\email{michael-david.chekroun@weizmann.ac.il}  
\email{mchekroun@atmos.ucla.edu}
\affiliation{Department of Atmospheric and Oceanic Sciences, University of California, Los Angeles, CA 90095-1565, USA}
\affiliation{Department of Earth and Planetary Sciences, Weizmann Institute of Science, Rehovot 76100, Israel} 

\author{Honghu Liu}
%\email{hhliu@vt.edu}
\affiliation{Department of Mathematics, Virginia Tech, Blacksburg, VA 24061, USA}

\author{James C. McWilliams}
\affiliation{Department of Atmospheric and Oceanic Sciences and Institute of Geophysics and Planetary Physics, University of California, Los Angeles, CA 90095-1565, USA}

\begin{abstract}%%%%%%%%%%%%
 A general, variational approach to derive low-order reduced models from possibly non-autonomous systems, is presented. 
The approach is based on the concept of optimal parameterizing manifold (OPM)
that substitutes the more classical notions of invariant or slow manifolds when breakdown of ``slaving'' occurs, i.e. when the unresolved variables cannot be expressed as an exact functional of the resolved ones anymore. The OPM provides, within a given class of  parameterizations of the unresolved variables, the manifold that averages out optimally these variables as conditioned on the resolved ones. 

The class of parameterizations retained here is that of continuous deformations of parameterizations rigorously valid near the onset of instability. These deformations  are  produced through integration of auxiliary backward-forward  (BF) systems built from the model's equations and lead to analytic formulas for parameterizations. In this modus operandi, the backward integration time is the key parameter to select per scale/variable to parameterize in order to derive the relevant  parameterizations which are doomed to be no longer exact away from  instability onset, due to breakdown of slaving typically encountered e.g.~for chaotic regimes. The selection criterion is then made through data-informed minimization of a least-square parameterization defect.  It is thus shown, through  optimization of the backward integration time per scale/variable to parameterize, that skilled OPM reduced systems can be derived for  predicting with accuracy higher-order critical transitions or catastrophic tipping phenomena, while training our parameterization formulas for regimes prior to these transitions take place.
\end{abstract}

%\date{\today}%
\date{September 15, 2023}%
\maketitle

{\bf 
We introduce a framework for model reduction to produce analytic formulas to parameterize the neglected variables. These parameterizations are built from the model's equations in which only a scalar is left to calibration per scale/variable to parameterize. This calibration is accomplished through a data-informed minimization of a least-square parameterization defect. By their analytic fabric the resulting parameterizations benefit physical interpretability.
Furthermore, our hybrid framework---analytic and data-informed--enables us to bypass the so-called extrapolation problem, known to be an important issue for purely data-driven machine-learned parameterizations. Here, by training our parameterizations prior transitions take place, we are able to perform accurate predictions of these transitions via the corresponding reduced systems. 
}

\section{Introduction}

This article is concerned with the efficient reduction of forced-dissipative systems of  the 
form 
\be\label{Eq_ODE_gen}
\frac{\d y}{\d t} = A y + G(y) +F(t), \qquad y \in H,
\ee
in which $H$ denotes a Hilbert state space, possibly of finite dimension. The forcing term $F$ is considered to be either constant in time or time-dependent, while $A$ denotes a linear operator not necessarily self-adjoint \footnote{but diagonalizable over the set of complex numbers $\mathbb{C}$}. that include dissipative effects, and $G$ denotes a nonlinear operator which saturate the possible unstable growth due to unstable modes and  may account for a loss of regularity (such as for nonlinear advection) in the infinite-dimensional setting. 
Such systems arise in broad range of applications; see e.g.~\cite{Hen81,GhCh87,Hale88,crawford1991introduction,CH93,SY02,MW05,KR11,temam2012infinite,MW14}.

The framework adopted is that of \citep{CLM19_closure} which allows for deriving analytic parameterizations of the neglected scales that represent efficiently the nonlinear interactions with the retained, resolved variables.  
The originality of the approach of \citep{CLM19_closure} is that the parameterizations are hybrid in their nature in the sense that they are both model-adaptive, based on the dynamical equations, and  data-informed by high-resolution simulations.

For a given system, the optimization of these parameterizations benefits thus from their analytical origin resulting  into only a few parameters to learn over short-in-time training intervals, mainly one scalar parameter per scale to parameterize. Their analytic fabric contributes then to their physical  interpretability compared to e.g.~parameterizations that would be machine learned. The approach is applicable to deterministic systems, finite- or infinite-dimensional, and is based on the concept of optimal parameterizing manifold (OPM) that substitutes the more classical notion of slow or invariant manifolds when takes place a breakdown of ``slaving'' relationships between the resolved and unresolved variables \citep{CLM19_closure},  i.e.~when the latter cannot be expressed as an exact functional of the formers anymore. 

By construction, the OPM takes its origin in a variational approach.
The OPM is the manifold that averages out optimally the neglected variables as conditioned on the the current state of the resolved ones \cite[Theorem 4]{CLM19_closure}.  The OPM allows for computing approximations  of the conditional expectation term arising in  the Mori-Zwanzig approach to stochastic modeling of neglected variables \cite{Chorin_al02,GKS04,gottwald2010recent,gottwald2017stochastic,LC23_review}; see also \cite[Theorem 5]{CLM19_closure} and \cite{chekroun2017emergence,chekroun2021stochastic}.

The approach introduced in \cite{CLM19_closure}  to determine OPMs in practice  consists at first deriving analytic parametric  formulas that match rigorous leading approximations of unstable/center manifolds or slow manifolds near e.g.~the onset of instability \cite[Sec.~2]{CLM19_closure}, and then to perform a data-driven optimization of the manifold formulas' parameters to handle regimes further away from that instability onset  \cite[Sec.~4]{CLM19_closure}.  
In other words, efficient parameterizations away from the onset are obtained as continuous deformations of parameterizations near the onset; deformations that are optimized by minimization of cost functionals tailored upon the dynamics and measuring the defect of parameterization.

There, the optimization stage allows for alleviating the small denominator problems rooted in small spectral gaps and for improving the parameterization of small-energy but dynamically important variables. Thereby, relevant parameterizations are derived in regimes where constraining spectral gap or timescale separation conditions are responsible for the well-known failure of standard invariant/inertial or slow manifolds \cite{CFNT89,debussche1991inertial,DebTem94,temam2010stability,temam2011slow,zelik2014inertial,CLM19_closure}.
As a result, the OPM approach provides (i) a natural remedy to  the excessive backscatter transfer of energy to the low modes classically encountered in turbulence  \cite[Sec.~6]{CLM19_closure}, and (ii)  allows for deriving optimal models of the slow motion for fast-slow systems not necessarily in presence of timescale separation \cite{chekroun2017emergence,chekroun2021stochastic}.  Due to their optimal nature, OPMs allow also for providing accurate parameterizations of dynamically important small-energy variables; a well-known issue encountered in the closure of chaotic dynamics and related to (i).

This work examine the ability of the theory-guided and data-informed parameterization approach of \cite{CLM19_closure} in deriving reduced systems able to predict higher-order transitions or catastrophic tipping phenomena, when the original, full system is possibly subject to time-dependent perturbations. From a data-driven perspective, this problem is tied to the so-called extrapolation problem, known for instance to be an  important issue in machine learning, requiring more advanced methods such as e.g.~transfer learning \cite{subel2022explaining}. While the past few decades have witnessed a surge and advances of many data-driven reduced-order modeling methodologies \cite{ahmed2021closures,hesthaven2022reduced}, the prediction of non-trivial dynamics for parameter regimes not included in the training dataset, remains a challenging task. Here, the OPM approach by its hybrid framework---analytic and data-informed---allows us to bypass this extrapolation problem at a minimal cost in terms of learning efforts as illustrated in Secns.~\ref{Sec_tippingmain} and ~\ref{Sec_RB}. As shown below, the training of OPMs at parameter values prior the transitions take place, is sufficient to perform accurate predictions of these transitions via  OPM reduced systems. 

The remainder of this paper is organized as follows. We first survey in Sec.~\ref{Sec_var} the OPM approach and provide the general variational parameterization formulas for  model reduction in presence of a time-dependent forcing. We then expand in Sec.~\ref{Sec_LIAmain} on the backward-forward (BF) systems method \cite{CLW15_vol2,CLM19_closure} to derive these formulas, clarifying fundamental relationships with homological equations arising in normal forms and invariant manifold theories \cite{cabre2005parameterization,haro2016parameterization,CLM19_closure,roberts2018backwards}. Section \ref{Sec_explicit} completes this analysis by analytic versions of these formulas in view of applications.  
The ability of predicting noise-induced transitions and catastrophic tipping phenomena \cite{kuehn2011mathematical,ashwin2012tipping} through OPM reduced systems is illustrated in Sec.~\ref{Sec_tippingmain} for a system arising in the study of thermohaline circulation. Successes in predicting higher-order transitions such as period-doubling and chaos are reported in Sec.~\ref{Sec_RB} for a Rayleigh-B\'enard problem, and contrasted by comparison with  the difficulties encountered by standard manifold parameterization approaches in Appendix \ref{Sec_IMfailue}. We then summarize the findings of this article with some concluding remarks in Sec.~\ref{Sec_conclusion}.

\section{Variational Parameterizations for Model Reduction}\label{Sec_var}
We summarize in this Section the notion of variational parameterizations introduced in \cite{CLM19_closure}.
The state space is decomposed as the sum of the subspace, $H_\c$, of resolved variables (``coarse-scale''),  and the subspace, $H_\s$, of unresolved variables (``small-scale''). In practice $H_\c$ is spanned by the first few eigenmodes of $A$ with dominant real parts (e.g.~unstable), and $H_\s$ by the rest of the modes, typically stable.  

In many applications, one is interested in deriving reliable reduced systems of Eq.~\eqref{Eq_ODE_gen} for wavenumbers smaller than a cutoff scale $k_c$ corresponding to a reduced state space $H_\c$ of dimension $m_c$. To do so, we are after parameterizations for the unresolved variables $y_\s$ (i.e., the vector component of the solution $y$ to  Eq.~\eqref{Eq_ODE_gen} in $H_\s$) of the form
\be\label{Eq_Htau}
\Phi_{\bftau}(X,t)= \sum_{n\geq m_c+1} \Phi_n(\tau_n,X,t) \boldsymbol{e}_n, \quad X \in H_\c, 
\ee
in which ${\bftau}$ the (possibly infinite) vector formed by the scalars $\tau_n$ is a free parameter. 
As it will be apparent below, the goal is to get a small-scale parameterization which is not necessarily exact such that $y_\s(t)$ is approximated by $\Phi_{\bftau}(y_\c(t),t)$ in a least-square sense, where $y_\c$ is the coarse-scale component of $y$ in $H_\c$. 
The vector $\boldsymbol{\tau}$ is aimed to be a homotopic deformation parameter.   The purpose is to have, as $\boldsymbol{\tau}$ is varied,  parameterizations that cover situations for which slaving relationships between the large and small scales hold ($y_\s(t)=\Phi_{\bftau}(y_\c(t),t))$ as well as situations in which they are broken ($y_\s(t)\neq \Phi_{\bftau}(y_\c(t),t))$, e.g.~far from criticality.
 With the suitable $\boldsymbol{\tau}$, the goal is to dispose of a reduced system resolving only the ``coarse-scales,'' able to e.g.~predict higher-order transitions caused by nonlinear interactions  with the neglected, ``small-scales.''

As strikingly shown in \cite{CLM19_closure}, meaningful parameterizations  away from the instability onset can still be obtained  by relying on bifurcation and center manifold theories. To accomplish this feat, the parameterizations, rigorously valid near that onset, need though to be revisited as suitably designed continuous deformations.
The modus operandi to provide such continuous deformations is detailed in Sec.~\ref{Sec_LIAmain} below, based on the method of backward-forward systems introduced in \cite[Chap.~4]{CLW15_vol2}, in a stochastic context. In the case where $G(y)=G_k(y) + h.o.t$ with $G_k$ a homogeneous polynomial of degree $k$), this approach gives analytical formulas of parameterizations given by (see Sec.~\ref{Sec_BF_gene})
%%%%%%%%%%%%%%%%%%%%%%%%%
\begin{widetext}
\bea\label{Eq_app2}
&\Phi_{\boldsymbol{\tau}}(X,t) =\sum_{n\geq m_c+1}\bigg( \int_{-\tau_n}^0e^{-s \lambda_n} \Big(\Pi_n G_k\big(
e^{s A_\c }X   - \eta(s) \big) + \Pi_n F (s+t) \Big) \,\mathrm{d}s \bigg) \boldsymbol{e}_n, \\
& \mbox{with }  X=\sum_{j=1}^{m_c} X_j \bm{e}_j\in H_{\c} \mbox{ and } \eta(s)= \int_s^0 e^{A_\c (s- r)} \Pi_c  F(r +t)  \d r.
\eea
\end{widetext}
Here, the $(\lambda_n,\boldsymbol{e}_n)$ denote the spectral elements of $A$ which are ordered such that $\Re(\lambda_n) \ge \Re(\lambda_{n+1})$.  In Eq.~\eqref{Eq_app2}, $\Pi_n$ and $\Pi_c$ denote the projector onto span$(\boldsymbol{e}_n)$ and $H_{\c}=$span$(\boldsymbol{e}_1,\cdots,\boldsymbol{e}_{m_c})$, respectively, while $A_\c=\Pi_\c A$. This formula provides an explicit expression for $\Phi^n_{\tau_n}(X,t)$ in \eqref{Eq_Htau}  given here by  the integral term over $[-\tau_n,0]$ multiplying $\boldsymbol{e}_{n}$ in \eqref{Eq_app2}. The only free parameter per mode $\bm{e}_n$ to parameterize is the backward integration time  $\tau_n$.

%%%%%%%%%%%%%%%%%%%%%%%%%
The vector ${\boldsymbol{\tau}}$, made of the $\tau_n$, is then optimized by using data from the full model. To allow for a better flexibility, the optimization is executed mode by mode one wishes to parameterize.
Thus, given a solution $y(t)$ of the full model Eq.~\eqref{Eq_ODE_gen} over an interval $I_T$ of length $T$, each parameterization $\Phi^n_{\tau_n}(X,t)$ of the solution amplitude $y_n(t)$ carried by mode  $\boldsymbol{e}_n$  is  optimized  by minimizing---in the $\tau_n$-variable---the following {\it parameterization defect}
\be\label{Eq_minQnHn}
\mathcal{Q}_n(\tau_n,T)= \overline{\big|y_n(t) -\Phi_n(\tau_n,y_\c(t),t) \big|^2},
\ee
for each $n\geq m_c +1$. Here $\overline{(\cdot )}$ denotes the time-mean over $I_T$ while $y_n(t)$ and $y_\c(t)$ denote the projections of $y(t)$ onto the high-mode $\boldsymbol{e}_n$ and the reduced state space $H_\c$, respectively.

%%%%%%%%%%%%%%%%%%%%
Geometrically,  the function $\Phi_{\bftau}(X,t)$ given by \eqref{Eq_app2} provides a time-dependent manifold $\mathfrak{M}_{\bftau}(t)$ such that:
\be\label{Eq_geom_interp0}
\overline{ \textrm{dist}(y(t),\mathfrak{M}_{\bftau}(t))^2} \leq \sum_{n=m_c+1}^N \mathcal{Q}_n(\tau_n,T),
\ee 
where $\textrm{dist}(y(t),\mathfrak{M}_{\bftau}(t))$ denotes the distance of $y(t)$ (lying e.g.~on the system's global attractor) to the manifold $\mathfrak{M}_{\bftau}$.

\begin{figure}%[hbtp]
\centering
\includegraphics[width=1\linewidth, height=.2\textwidth]{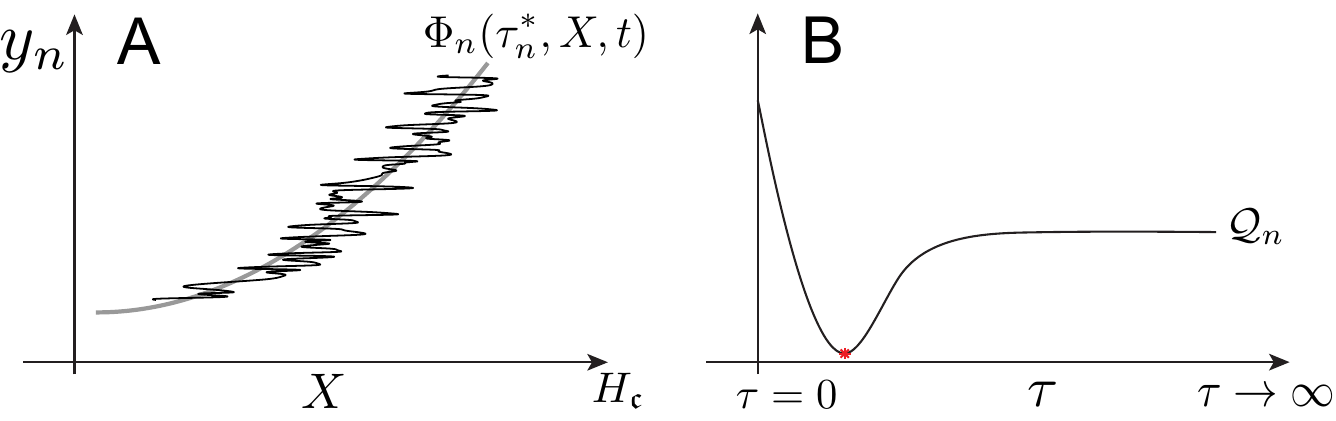}
%\vspace{-1ex}
\caption{{\bf Panel A:}  The black curve represents the training trajectory, here shown to to be transverse to the parameterizing manifolds (i.e.~absence of exact slaving). The grey smooth curves represent the time-dependent OPM aimed at tracking the state of the neglected variable $y_n$ at time $t$ as a function of the resolved variables $X$.  {\bf Panel B:}  A schematic of the dependence on $\tau$ of the parameterization defect $\mathcal{Q}_n$ given by \eqref{Eq_minQnHn}. The red asterisk marks the minimum of  $\mathcal{Q}_n$ achieved for $\tau=\tau_n^\ast$.}\label{Fig_intro}
\end{figure}

Thus minimizing each $\mathcal{Q}_n(\tau_n,T)$ (in the $\tau_n$-variable) is a natural idea to enforce closeness of $y(t)$ in a least-squares sense to  the manifold $\mathfrak{M}_{\bftau}(t)$.
Panel A in Fig.~\ref{Fig_intro} illustrates \eqref{Eq_geom_interp0} for the $y_n$-component: The  optimal parameterization, $\Phi_n(\tau_n^\ast,X,t)$, minimizing \eqref{Eq_minQnHn} is shown for a case where the dynamics is transverse to it (e.g.~in absence of slaving) while $\Phi_n(\tau_n^\ast,X,t)$ provides the best parameterization in a least-squares sense.

%%%%%%%%%%%%%%%%%%%%%%%%%%%%%
%%%%%%%%%%%%%%%%%%%%%%%%%%%%%

In practice, the {\it normalized parameterization defect}, $Q_n$, is often used to compare different parameterizations.  
It is defined as $Q_n(\tau,T)=\mathcal{Q}_n(\tau,T)/\overline{|y_n|^2}$.  For instance, the flat manifold corresponding to  no parameterization ($\tau_n=0$) of the neglected variables (Galerkin approximations) comes with    $Q_n=1$ for all $n$, while a manifold corresponding to a perfect slaving relationship between the $y_\c(t)$ and  $y_n(t)$'s, comes with  $Q_n=0$. 
 When $0<Q_n<1$ for all $n$, the underlying manifold $\mathfrak{M}_{\bftau}$ will be referred to as a {\it parameterizing manifold (PM)}.  
 Once the parameters $\tau_n$ are optimized by minimization of \eqref{Eq_minQnHn}, the resulting manifold will be referred to as the the optimal parameterizing manifold (OPM) \footnote{ Here, we make the abuse of language that the OPM is sought within the class of parameterizations obtained through the backward-forward systems of Sec.~\ref{Sec_LIAmain}. In general, OPMs are more general objects related to the conditional expectation associated with the projector $\Pi_\c$; see \cite[Sec.~3.2]{CLM19_closure}.}.
 
We conclude this section by a few words of practical considerations. As documented in \cite[Secns.~5 and 6]{CLM19_closure}, the amount of training data required in order to  reach a robust  estimation of the optimal backward integration time $\tau^*_n$,   is often comparable to the dynamics' time horizon that is necessary to resolve the decay of correlations for the high-mode amplitude $y_n(t)$. For multiscale systems, one thus often needs  to dispose of a training dataset sufficiently large to resolve the slowest decay of temporal correlations of the scales to parameterize. On the other hand, by benefiting from their dynamical origin, i.e.~through the model's equations, the parameterizations formulas  employed in this study (see Sec.~\ref{Sec_LIAmain})  allow often  for reaching out in practice satisfactory OPM approximations when optimized over training intervals shorter than these characteristic timescales.
   
When these conditions are met, the minimization of the parameterization defect \eqref{Eq_minQnHn} is performed 
by a simple gradient-descent algorithm \cite[Appendix]{CLM19_closure}. There, the first local minimum that one reaches, 
corresponds often to the OPM; see Secns.~\ref{Sec_tippingmain} and ~\ref{Sec_RB} below, and \cite[Secns.~5 and 6]{CLM19_closure}.  
In the rare occasions where the parameterization defect exhibits more than one local minimum and the corresponding local minima are close to each others, criteria involving colinearity between the features to parameterize  and the parameterization itself can be designed to further assist the OPM selection. Section \ref{Sec_localmin} below illustrates this point with the notion of parameterization correlation.

%%%%%%%%%%%%%%%%%%%%%%%%%%%%%%%%%%%%
\section{Variational Parameterizations: Explicit Formulas}\label{Sec_LIAmain}

\subsection{Homological Equation and Backward-Forward Systems: Time-dependent Forcing}\label{Sec_BF_gene}

In this section, we recall from \cite{CLM19_closure}  and generalize to the case of time-dependent forcing, the theoretical underpinnings of the parameterization formula \eqref{Eq_app2}. 
First, observe that when $F \equiv 0$ and $\tau_n = \infty$ for all $n$, the parameterization \eqref{Eq_app2} is reduced (formally) to 
\be \label{Eq_LP_integral}
\mathfrak{J}(X) = \int_{-\infty}^0 e^{-sA_{\s}} \Pi_{\s} G_k(e^{sA_{\c}}X) \d s, \; X\in H_\c,
\ee
where  $A_{\s} = \Pi_{\s} A$, with $\Pi_{\s}$ denoting the canonical projectors onto  $H_{\s}$. Eq.~\eqref{Eq_LP_integral} is known in  invariant manifold theory as a Lyapunov-Perron integral \cite{SY02}. It provides the leading-order approximation of the underlying invariant manifold function if a suitable spectral gap condition is satisfied
and solves the homological equation:
\be\label{Fundamental_Eq}
D\psi (X)A_\c X-A_\s \psi(X)=\Pi_\s G_k(X),
\ee
see \cite[Lemma 6.2.4]{Hen81} and \cite[Theorem 1]{CLM19_closure}. 
The later equation is a simplification of the invariance equation providing the invariant manifold when it exists; see \cite[Sec.~VIIA1]{crawford1991introduction} and  \cite[Sec.~2.2]{CLM19_closure}.
Thus, Lyapunov-Perron integrals such as \eqref{Eq_LP_integral} are intimately related to the homological equation, and the study of the latter informs us on the former, and  in turn on the closed form of the underlying invariant manifold. 

Another key element was pointed out  in \cite[Chap.~4]{CLW15_vol2}, in the quest of getting new insights about invariant manifolds in general and more specifically concerned with the approximation of stochastic invariant manifolds of stochastic PDEs \cite{CLW15_vol1}, along with the rigorous low-dimensional approximation of their dynamics \cite{chekroun2023transitions}. It consists of the method of Backward-Forward (BF) systems, thereafter revealed  in  \cite{CLM19_closure} as a key tool to produce parameterizations (based on model's equations) that are relevant beyond the domain of applicability of invariant manifold theory, i.e.~away the onset of instability.

To better understand this latter feature, let us recall first that BF systems allow for providing to Lyapunov-Perron integrals such as \eqref{Eq_LP_integral}, a flow interpretation. 
In the case of \eqref{Eq_LP_integral}, this BF system takes the form \cite[Eq.~(2.29)]{CLM19_closure}:
\begin{subequations} \label{Eq_BF}
\begin{align}
& \frac{\mathrm{d} p}{\d s} =  A_{\c} p, \hspace{6.7em} s \in[ -\tau, 0],    \label{BF1} \\
& \frac{\mathrm{d} q}{\d s} = A_{\s} q  +  \Pi_{\s} G_k(p),  \qquad  s \in [-\tau, 0], \label{BF2}\\
& \mbox{with } p(s)\vert_{s=0} = X, \mbox{ and } q(s)\vert_{s=-\tau}=0.
\end{align}
\end{subequations}
The Lyapunov-Perron integral is indeed recovered from this BF system. 
The solution to Eq.~\eqref{BF2} at $s=0$ is given by
\be\label{Eq_h1_tau}
q(0) =\int_{-\tau}^0 e^{-sA_{\s}} \Pi_{\s} G_k(e^{sA_{\c}}X) \d s, \; X\in H_\c,
\ee
which by taking the limit formally in \eqref{Eq_h1_tau} as $\tau\rightarrow \infty$, leads to $\mathfrak{J}$ given by \eqref{Eq_LP_integral}.  Thus, stretching $\tau$ to infinity in the course of integration of the BF system \eqref{Eq_BF}  allows for recovering rigorous approximations (under some non-resonance conditions  \cite[Theorem 1]{CLM19_closure}) of well-known objects such as the center manifold; see also \cite[Lemma 6.2.4]{Hen81} and \cite[Appendix A.1]{MW14}. 
The intimate relationships between Lyapunov-Perron integral, $\mathfrak{J}(X)$, and the homological Eq.~\eqref{Fundamental_Eq} are hence supplemented by their relationships with the BF system  \eqref{Eq_BF}.

By breaking down, mode by mode, the backward integration of Eq.~\eqref{BF2} (in the case of $A_\s$ diagonal), one allows for a backward integration time $\tau_n$ per mode's amplitude to parameterize, leading in turn to 
 the following class of parameterizations 
\be\label{Eq_app0}
\Phi_{\boldsymbol{\tau}}(X) =\hspace{-1ex}\sum_{n\geq m_c+1} \hspace{-1ex} \bigg( \int_{-\tau_n}^0e^{-s \lambda_n} \Pi_n G_k\big(
e^{s A_\c}X  \big) \,\mathrm{d}s \bigg) \boldsymbol{e}_n, 
\ee
as indexed by $\boldsymbol{\tau}=(\tau_n)_{n\geq m_c+1}$. Clearly Eq.~\eqref{Eq_app2} is a generalization of Eq.~\eqref{Eq_app0}.

 We make precise now the similarities and differences  between Eq.~\eqref{Eq_app2} and Eq.~\eqref{Eq_app0}, at a deeper level as informed by their BF system representation. In that respect, we consider for the case of time-dependent forcing, the following BF system 
 \begin{subequations} \label{Eq_BF_withforcing}
\begin{align}
& \frac{\mathrm{d} p}{\d s} =  A_{\c} p + \Pi_{\c}{F}(s), \hspace{5.5em} s  \in[ t-\tau, t],    \label{BF1_withforcing} \\
& \frac{\mathrm{d} q}{\d s} = A_{\s} q  +  \Pi_{\s} G_k(p) +\Pi_{\s}{F}(s),  \;\;  s \in[ t-\tau, t], \label{BF2_withforcing}\\
& \mbox{with } p(t) = X \in H_\c, \;  \mbox{ and } q(t-\tau)=0. 
\end{align}
\end{subequations}
Note that compared to the BF system \eqref{Eq_BF}, the backward-forward integration is operated here on $[ t-\tau, t]$  to account for the non-autonomous character of the forcing, making this way the corresponding parameterization time-dependent as summarized in Lemma \ref{Lemma_foundation} below. In the presence of an autonomous forcing, operating the backward-forward integration of the  BF system \eqref{Eq_BF} on $[ t-\tau, t]$ does not change the parameterization, which remains time-independent.

Since the BF system \eqref{Eq_BF_withforcing} is one-way coupled (with $p$ forcing \eqref{BF2_withforcing} but not reciprocally), if one assumes $A_{\s}$ to be diagonal over $\mathbb{C}$, one can break down the forward equation Eq.~\eqref{BF2_withforcing} into the  equations,
\be\label{qn_eq}
 \frac{\mathrm{d} q_n}{\d s} = \lambda_n q_n  +  \Pi_{n} G_k(p) +  \Pi_{n}{F}(s),  \;\;  s  \in[ t-\tau, t], 
\ee
allowing for flexibility in the choice of $\tau$ per mode $\boldsymbol{e}_n$ whose amplitude is aimed at being  parameterized by $q_n$, for each $n\geq m_c+1$. 

After backward integration  of Eq.~\eqref{BF1_withforcing} providing $p(s)$, one obtains by forward integration of Eq.~\eqref{qn_eq} that 
\bea \label{eq_h1_NDS0}
 q_n(t)=\int_{t-\tau}^t e^{(t-s')\lambda_n} &\Pi_{n} G_k\big(e^{(s'-t) A_{\c}} X-\gamma(s')\big) \d s' \\
 &+\int_{t-\tau}^te^{(t-s') \lambda_n} \Pi_n F(s')\d s',
 \eea
 with $\gamma(s')=\int_{s'}^t e^{(s'-r) A_{\c}}\Pi_\c F(r)\d r$.  By making the change of variable $s'=s+t$, one gets
\bea  \label{Eq_qn_withF}
 q_n(t)= \int_{-\tau}^0 e^{-s\lambda_n} & \Pi_{n} G_k\Big(e^{s A_{\c}}X -\eta(s)\Big) \d s\\
  & + \int_{-\tau}^0 e^{-s \lambda_n} \Pi_n F(s+t)\d s,
\eea
with $\eta(s)=\int_s^0 e^{(s-r)A_\c}\Pi_\c F(r+t) \d r =\gamma(s+t)$. Now by summing up these parameterizations over the modes $\boldsymbol{e}_n$ for $n\geq m_c+1$, we arrive at 
Eq.~\eqref{Eq_app2}. Thus, our general parameterization formula \eqref{Eq_app2} is derived by solving the BF system made of backward Eq.~\eqref{BF1_withforcing}  and the forward Eqns.~\eqref{qn_eq}.

%%%%%%%%%%%%%%
We want to gain into interpretability of such parameterizations in the context of forced-dissipative systems such as Eq.~\eqref{Eq_ODE_gen}.
For this purpose, we clarify the fundamental equation satisfied by our parameterizations; the goal being here to characterize 
the analogue of \eqref{Fundamental_Eq} in this non-autonomous context. 
To simplify, we restrict to the case $\Pi_\c F=0$.  The next Lemma provides the sought equation. 
%The next Lemma provides useful elements for interpretability. 

%%%%%%%%%%%%%%
\begin{widetext}
\begin{lem}\label{Lemma_foundation} 
Assume $\Pi_\c F=0$ in the BF system  \eqref{Eq_BF_withforcing}. The solution $q_{X,\tau}(t)$ to Eq.~\eqref{BF2_withforcing} is given by
\be \label{eq_h1_NDS}
q_{X,\tau}(t)= \int_{-\tau}^0 e^{-s A_{\s}} \Pi_{\s} G_k(e^{s A_{\c}}X) \d s + \int_{-\tau}^0 e^{-s A_{\s}} \Pi_\s F(s+t)\d s.
\ee
It provides a time-dependent manifold function that  maps $H_\c$ into $H_\s$ and satisfies the following  first order system of PDEs:
\be \label{Eq_invariance_NDS}
 \Big(\partial_t+ \mathcal{L}_A\Big) \Phi (X,t)=\underbrace{\Pi_{\s} G_k(X)- e^{\tau A_{\s}} \Pi_{\s} G_k(e^{-\tau A_{\c}}X)}_{(I)} +\underbrace{\Pi_\s F(t)- e^{\tau A_{\s}}  \Pi_{\s} F(t-\tau)}_{(II)},
\ee
with
$\mathcal{L}_A$ being the differential operator acting on smooth mappings $\psi$ from $H_\c$ into $H_\s$, defined as:
\be \label{Def_L}%\label{h1_eqn}
\mathcal{L}_A[\psi] (X)=D \psi(X) A_{\c} X  - A_{\s} \psi(X), \; X\in H_\c.
\ee
\end{lem}
The proof of this Lemma is found in Appendix \ref{Proof_Lemma}. 

\end{widetext}
%%%%%%%%%%%%%%%%%
As it will become apparent below, the $\tau$-dependence of the terms in (I) is meant to control the small denominators that arise in presence of small spectral gap between the spectrum of $A_\c$ and $A_\s$ that leads typically to over-parameterization when standard invariant/inertial manifold theory is applied in such situations; see Remark \ref{Rem_smallgap} below and \cite[Sec.~6]{CLM19_closure}.

The terms in (II) are responsible for the presence of exogenous memory terms in the solution to the homological equation Eq.~\eqref{Eq_invariance_NDS}, i.e.~in the parameterization $\Phi(X,t)$; see Eq.~\eqref{Eq_param_timedependent} below.

In the case $F(t)$ is bounded \footnote{One could also express a constraint on the behavior of $F(t)$ as $t \rightarrow \infty$. For instance one may require that  
\be
\lim_{s\rightarrow -\infty}e^{A_{\c}s} \int_{s}^0 e^{A_{\c}s'} \Pi_{\c} F (s') \d s'<\infty, etc.
\ee}, 
and $\Re \sigma (A_\s)<0$, Eq.~\eqref{Eq_invariance_NDS} reduces, in the (formal) limit $\tau \rightarrow \infty$, to:
\be\label{Eq_homological_NDS}
 \Big(\partial_t+ \mathcal{L}_A\Big) \Phi (X,t)=G_k(X)+\Pi_\s F(t).
\ee
Such a system of linear PDEs is known as the homological equation and arises  in the theory of time-dependent normal forms \cite{Arnold88,haro2016parameterization}; see also \cite{PR06}. 
In the case where $\Pi_\s F$ is $T$-periodic, one can seek an approximate solution in a Fourier expansion form  \cite[Sec.~5.2]{sinha2006bifurcation} leading to useful insights. For instance, solutions to \eqref{Eq_homological_NDS} when $G$ is quadratic, exist if the following non-resonance  condition is satisfied
 \bea
& \mathfrak{i} \nu \frac{\pi}{T}+\lambda_{i}+\lambda_{j}-\lambda_n \neq 0, \; \nu \in \mathbb{Z},\\
&\text{for } (i,j)\in I^2, \; n\geq m_c+1,\; 
 \eea
 (with $I=\{1,\cdots,m_c\}$ and $\mathfrak{i}^2=-1$)
 in the case $\sigma(A)$ is discrete, without Jordan block.  
 Actually one can even prove in this case that
 \be
\mbox{Spec}(\partial_t+ \mathcal{L})=\bigg\{ \mathfrak{i} \nu \frac{\pi}{T}+\delta_{ij}^n,\; (i,j)\in I^2, \; n\geq m_c+1,\; \nu \in \mathbb{Z} \bigg\},
 \ee
 (where $\delta_{ij}^n=\lambda_{i} + \lambda_{j} - \lambda_{n}$)
 on the space of functions 
 \be
\hspace{-1.5ex} \mathcal{E}=\left\{ \hspace{-.5ex}\sum_{n \geq m_c+1} \hspace{-.8ex} \bigg( \sum_{i,j=1}^{m_c} \Gamma^n_{ij}(t) X_i X_j \bigg) \boldsymbol{e}_n, \; \Gamma^n_{i j} \in {L}^2(\mathbb{T}) \right\},
 \ee
in which $X_i$ and $X_j$ represent the components of the vector $X$ in $H_\c$ onto $\boldsymbol{e}_i$ and  $\boldsymbol{e}_j$, respectively. 

Thus, in view of Lemma \ref{Lemma_foundation} and what precedes, the small-scale parameterizations \eqref{eq_h1_NDS} obtained by solving the BF system  \eqref{Eq_BF_withforcing} over finite time intervals, can be  conceptualized as perturbed solutions to the homological equation Eq.~\eqref{Eq_homological_NDS} arising in the computation of invariant and normal forms of non-autonomous dynamical systems \cite[Sec.~2.2]{haro2016parameterization}. 
The perturbative terms brought by the $\tau$-dependence play an essential role to cover situations beyond the domain of application of normal form and invariant manifold theories. As explained in Sec.~\ref{Sec_explicit} and illustrated in Secns.~\ref{Sec_tippingmain}-\ref{Sec_RB} below, these terms can be optimized to ensure skillful parameterizations for predicting,  by reduced systems, higher-order transitions escaping the domain of validity of these theories.

\subsection{Homological Equation and Backward-Forward Systems: Constant Forcing}\label{Sec_BF_gene2}
 In this Section, we clarify the conditions under which solutions to Eq.~\eqref{Eq_BF_withforcing} exist in the asymptotic limit $\tau \rightarrow \infty$. We consider the case where $F$  is constant in time, to simplify. For the existence of Lyapunov-Perron integrals under the presence of more general time-dependent forcings we refer to \cite{PR06}.
To this end, we  introduce, for $X$ in $H_\c$,
\be\label{Eq_h1_withforcing}
\mathfrak{J}_{F}(X) = \int_{-\infty}^0 e^{-sA_{\s}} \left( \Pi_{\s} G_k\big(p(s)\big) + \Pi_{\s}{F} \right) \d s,
\ee
where $p(s)$ is the solution to  Eq.~\eqref{BF1_withforcing}, namely 
\be \label{Eq_p_soln}
p(s) = e^{sA_{\c}}X -  \int_s^0 e^{A_{\c}(s- s')} \Pi_{\c} F \d s'.
\ee
We have then the following result that we cast for the case of finite-dimensional ODE system (of dimension $N$) to simply. 

\bt \label{Thm_BF}
Assume that $F$ is constant in time and that $A$ is diagonal under its eigenbasis $\{\bm{e_j} \in \mathbb{C}^N \;  | \;  j=1,\cdots N\}$. 
Assume furthermore that $\Re(\lambda_{m_c + 1}) < \Re(\lambda_{m_c})$, and $\Re \sigma(A_\s)<0$,  i.e.~$H_\s$ contains only stable eigenmodes. 

Finally, assume that the following  non-resonance condition holds:
\bea \label{Eq_NR}
&\forall \, (j_1,\cdots, j_k)\in (1,\cdots,m_c)^k, \; n \ge m_c+1,  \\
&\left( G_{j_1 \cdots j_k}^n \neq 0 \right)  \implies \left( \Re(\lambda_n-\sum_{p = 1}^k \lambda_{j_p}) < 0\right),
\eea
with $ G_{j_1 \cdots j_k}^n =\big\langle G_k(\bm{e}_{j_1},\cdots,\bm{e}_{j_k}), \bm{e}_n^\ast \big\rangle $, where $\langle \cdot, \cdot\rangle$ denotes the inner product on $\mathbb{C}^N$, $G_k$ denotes the leading-order term in the Taylor expansion of $G$, and $\bm{e}_n^\ast$ denotes the eigenvectors of the conjugate transpose of $A$.

 Then, the Lyapunov-Perron integral  $\mathfrak{J}_{F}$ given by \eqref{Eq_h1_withforcing} is well defined and is a solution to the following homological equation
 \be\label{Eq_homoligical}
\mathcal{L}_A[\psi] (X) + D \psi[X] \Pi_{\c} F = \Pi_{\s}G_k(X) + \Pi_{\s} F,
\ee
and provides the leading-order approximation of the invariant manifold function $h(X)$ in the sense that 
\bes
\norm{\mathfrak{J}_{F}(X)-h(X)}_{H_\s}=o(\norm{X}_{H_\c}^k),\; \; X \in H_\c.
\ees
Moreover, $\mathfrak{J}_{F}$ given by \eqref{Eq_h1_withforcing} is the limit, as $\tau$ goes to infinity, of the solution to the BF system \eqref{Eq_BF_withforcing}, when $F$ is constant in time. That is
\be \label{Eq_PB_limit}
\lim_{\tau \rightarrow \infty} \| q_{X,\tau}(0) - \mathfrak{J}_F(X) \| = 0,
\ee
where $q_{X,\tau}(0)$ is the solution to  Eq.~\eqref{BF2_withforcing}.
\et
%%%%%%%%%

Conditions similar to \eqref{Eq_NR} arise in the smooth linearization of dynamical systems near an equilibrium  \cite{sell1985smooth}. Here, condition \eqref{Eq_NR} implies that the eigenvalues of the stable part satisfy a Sternberg condition of order $k$ \cite{sell1985smooth} with respect to the eigenvalues associated with the modes spanning the reduced state space $H_\c$.

 This theorem is essentially a consequence of  \cite[Theorems 1 and 2]{CLM19_closure}, in which condition \eqref{Eq_NR} is a stronger version of that  used for \cite[Theorem 2]{CLM19_closure}; see also \cite[Remark 1 (iv)]{CLM19_closure}. This condition is necessary and sufficient here for $\int_{-\infty}^0 e^{-sA_{\s}} \Pi_{\s} G_k(p(s)) \d s$ to be well defined. The convergence of $\int_{-\infty}^0 e^{-sA_{\s}} \Pi_{\s}{F} \d s$ is straightforward since $\Re \sigma(A_\s)<0$ by assumption and $F$ is constant. The derivation of \eqref{Eq_homoligical} follows the same lines as the derivation for \cite[Eq.~(4.6)]{CLM19_closure}.

One can also generalize Theorem~\ref{Thm_BF} to the case when $F$ is time-dependent provided that $F$ satisfies suitable conditions to ensure $\mathfrak{J}_{F}$ given by \eqref{Eq_h1_withforcing} to be well defined and that the non-resonance condition \eqref{Eq_NR} is suitably augmented. We leave the precise statement of such a generalization to a future work. For the applications considered in later sections, the forcing term $F$ is either a constant or with a sublinear growth. For such cases, $\mathfrak{J}_{F}$ is always well defined under the assumptions of Theorem~\ref{Thm_BF}. We turn now to present the explicit formulas of the parameterizations based on the BF system \eqref{Eq_BF_withforcing}.

\subsection{Explicit Formulas for Variational Parameterizations}\label{Sec_explicit}
We provide in this Section, closed form formulas of parameterizations 
for forced dissipative such as Eq.~\eqref{Eq_ODE_gen}. We consider first the case where $F$ is constant in time, and then deal with time-dependent forcing case. 

\subsubsection{The constant forcing case}
To favor flexibility in applications,  we consider scale-awareness of our parameterizations via BF systems, i.e.~we consider parameterizations that are built up, mode by mode, from the integration, for each $n\geq m_c+1$, of the following BF systems
\begin{subequations} \label{Eq_BF_quad}
\begin{align}
& \frac{\mathrm{d} p}{\d s} =  A_\c  p + \Pi_{\c} F,     \hspace{6.7em} s \in[ -\tau, 0],\label{BF1_quad} \\
& \frac{\mathrm{d} q_{n}}{\d s} = \lambda_n q_n  +  \Pi_{n} G_k\big(p(s)\big) + \Pi_{n} F, \;  s \in [-\tau, 0] \label{BF2_quad}\\
& \mbox{with } p(0) = X \in H_\c, \mbox{ and } q_n(-\tau)=0.\label{BF2_quad3}
\end{align}
\end{subequations}
Recall that $q_n(0)$ is aimed at parameterizing the solution amplitude $y_n(t)$ carried  by mode  $\boldsymbol{e}_n$,  when $X=y_\c(t)$ in  Eq.~\eqref{BF2_quad3}, whose parameterization defect is minimized by optimization of the backward integration time $\tau$ in  \eqref{Eq_minQnHn}. To dispose of explicit formulas for $q_n(0)$ qualifying the dependence on $\tau$, facilitates greatly this minimization.  

In the case the nonlinearity $G(y)$ is quadratic, denoted by $B(y,y)$, such formulas are easily accessible as \eqref{Eq_BF_quad} is integrable. It is actually integrable in the presence of higher-order terms, but we provide the details here for the quadratic case, leaving the details to the interested reader for extension.

Recall that we denote by $\boldsymbol{e}_j^\ast$ the conjugate modes from the adjoint  $A^\ast$ and that these modes satisfy the bi-orthogonality condition. That is $\langle \boldsymbol{e}_i, \boldsymbol{e}_j^\ast\rangle = 1$ if $i = j$, and zero otherwise, where $\langle \cdot, \cdot \rangle$ denotes the inner product endowing $H$. 
Denoting by $\Phi_n(\tau, \boldsymbol{\lambda}, X)$ the solution $q_n(0)$ to Eq.~\eqref{BF2_quad}, we find after calculations that
\bea \label{Eq_Phi_tau}
\Phi_n(\tau, \boldsymbol{\lambda}, X) =R_n&(F,\bflambda,\tau,X) - \frac{1 - e^{\tau \lambda_n}}{\lambda_n} F_n, \\
&+ \sum_{i, j = 1}^{m_c} D_{ij}^n(\tau, \boldsymbol{\lambda}) B_{ij}^n X_{i} X_{j},
\eea
in which $X_i$ and $X_j$ denote the components of the vector $X$ in $H_\c$ onto $\boldsymbol{e}_i$ and  $\boldsymbol{e}_j$, $B_{ij}^n=\langle B(\boldsymbol{e}_i, \boldsymbol{e}_j), \boldsymbol{e}^*_n\rangle$ and  
\be \label{Eq_D_term0}
D_{ij}^n(\tau,\bflambda)= \frac{1 - \exp\big(-\delta_{ij}^n\tau\big)}{\delta_{ij}^n},  \text{if $\delta_{ij}^n\neq 0$},
\ee
while $D_{ij}^n(\tau,\bflambda)=\tau$ otherwise. Here $\delta_{ij}^n=\lambda_{i} + \lambda_{j} - \lambda_{n}$, with the $\lambda_j$ referring to the eigenvalues of $A$. 
We refer to \cite{CLM19_closure} and Appendix \ref{Sec_LIA_formula} for the expression of $R_n(F,\bflambda,\tau,X)$ which accounts for the nonlinear interactions between the forcing components in $H_\c$.  Here,  the dependence on the $\lambda_j$ in \eqref{Eq_Phi_tau} is made apparent, as this dependence plays a key role in the prediction of the higher-order transitions; see applications to the Rayleigh-B\'enard system of Sec.~\ref{Sec_RB}. 

%%%%%%%%%%%%%%%%%%%%%%
\br\label{Rem_smallgap}
{\bf OPM balances small spectral gaps.}
Theorem~\ref{Thm_BF} teaches us that when $\delta_{ij}^n>0$ not only $D_{ij}^n(\tau,\bflambda)$ defined in \eqref{Eq_D_term0} converges towards a well defined quantity as $\tau_n\rightarrow \infty$  but also the coefficients involved in  $R_n(F,\bflambda,\tau,X)$ (see \eqref{Formula_GammaF}-\eqref{Eq_Un}-\eqref{Eq_Vn}), in case of existence of invariant manifold. 
For parameter regimes where the latter fails to exist, some of the $\delta_{ij}^n$ or the $\lambda_j-\lambda_n$ involved in \eqref{Eq_Un}-\eqref{Eq_Vn} can become small, leading to the well known small spectral gap issue \cite{zelik2014inertial} manifested typically by over-parameterization of the $\bm{e}_n$'s mode amplitude when the Lyapunov-Perron parameterization \eqref{Eq_h1_withforcing} is employed; see \cite[Sec.~6]{CLM19_closure}. The presence of the $\tau_n$ through the exponential terms in \eqref{Eq_D_term0} and \eqref{Eq_Un}-\eqref{Eq_Vn} allows for balancing these small spectral gaps after optimization and improve notably the parameterization and closure skills; see Sec.~\ref{Sec_RB} and Appendix \ref{Sec_IMfailue}.  
\er

\br\label{Rem_FMT}
Formulas such as introduced in \cite{FMT88}  in approximate inertial manifold (AIM) theory, and used earlier in atmospheric dynamics \cite{daley1980development} in the context of non-normal modes initialization \cite{baer1977complete,machenhauer1977dynamics,leith1980nonlinear}, are also tied to leading-order approximations of invariant manifolds  since $\mathfrak{J}(X)$ given by  \eqref{Eq_LP_integral} satisfies
\be\label{Eq_FMTgen}
\mathfrak{J}(X) =-A_\s^{-1} \Pi_\s G_k(X) +\mathcal{O}(\norm{X}^k), \; X \in H_\c,
\ee
and the term $-A_\s^{-1} \Pi_\s G_k(X)$ is the main parameterization used in \cite{FMT88}; see \cite[Lemma 4.1]{CLW15_vol2} and \cite[Theorem A.1.1]{MW14}. 
Adding higher-order terms can potentially lead to more accurate parameterizations and push the validity of the approximation to a larger neighborhood \cite{roberts1988application,brown1991minimal}, but in presence of small spectral gaps, such an operation may  become also less successful.  
When such formulas arising in AIM theory are inserted within the proper data-informed variational approach such as done in \cite[Sec.~4.4]{CLM19_closure}, their  optimized version  can also handle regimes in which spectral gap becomes small as demonstrated in  \cite[Sec.~6]{CLM19_closure}.
\er

%%%%%%%%%%%%%%%%%%%%%%%%%
\subsubsection{The time-dependent forcing case}
Here, we assume that the neglected modes are subject to time-dependent forcing according to $F(t)=\sum_{n\geq m_c+1} \sigma_n f_n(t) \bm{e}_n$. Then by solving  the BF systems made of Eq.~\eqref{BF1_withforcing} and Eqns.~\eqref{qn_eq} with $q_n(t-\tau)=\zeta$ (to account for possibly non-zero mean), we arrive at the following time-dependent parameterization of the $n$th  mode's amplitude:
\begin{widetext}
\be\label{Eq_param_timedependent}
\Phi_n(\tau, \bm{\lambda}, X,t)=e^{\lambda_n \tau} \zeta+\underbrace{\sigma_n e^{\lambda_n t} \int_{t-\tau}^t e^{-\lambda_n s} f_n(s) \d s}_{(I)} + \sum_{i, j = 1}^{m_c} D_{ij}^n(\tau, \boldsymbol{\lambda}) B_{ij}^n X_{i} X_{j},
\ee
where $X=\sum_{j=1}^{m_c} X_j \bm{e}_j$ lies in $H_\c$.
\end{widetext}

This formula of $\Phi_n$  gives the solution to the homological equation Eq.~\eqref{Eq_invariance_NDS} of Lemma \ref{Lemma_foundation} in which $\Pi_\s$ therein is replaced here by $\Pi_n$, the projector onto the mode $\bm{e}_n$ whose amplitude is parameterized by $\Phi_n$. Clearly, the terms in (I) are produced by the time-dependent forcing.  They are functionals of the past $f_n$ and convey thus  {\it exogenous} memory effects. 

The integral  term in (I) of Eq.~ \eqref{Eq_param_timedependent} is of the form 
\be \label{Eq_def_I}
I(t)= e^{\kappa t} \int_{t-\tau}^t e^{-\kappa s} f_n(s) \d s.
\ee
By taking derivates on both sides of \eqref{Eq_def_I}, we observe that $I$ satisfies 
\be \label{Eq_for_I}  
\frac{\d I}{\d t} = \kappa I + f_n(t)- e^{\kappa \tau} f_{n}(t-\tau).
\ee
As a practical consequence, the computation of $I(t)$ boils down to solving the scalar ODE \eqref{Eq_for_I}  which can be done with high-accuracy numerically, when $\kappa < 0$. 
One bypasses thus the computation of many quadratures (as $t$ evolves) that we would have to perform when relying on Eq.~\eqref{Eq_def_I}. Instead only one quadrature is required corresponding to the initialization of  the ODE \eqref{Eq_for_I} at $t=0$.

This latter computational aspect is important not only for simulation purposes when the corresponding OPM reduced system is ran online  but also for training purposes, in the search of the optimal $\tau$ during the offline minimization stage of the parameterization defect.  

If time-dependent forcing terms are present in the reduced state space $H_\c$, then the BF system \eqref{Eq_BF_quad} can still be solved analytically albeit leading to more involved integral terms than in (I) in the corresponding parameterization. 
This aspect will be detailed elsewhere.

\subsubsection{OPM reduced systems}
Either in the constant forcing or time-dependent forcing case, our OPM reduced system 
takes the form 
\be\label{Eq_reduced}
\dot{X}=A_\c X +\Pi_\c G\big(X +\Phi_{\bm{\tau}^\ast}(X,t)\big)+\Pi_\c F,
\ee
where 
\be\label{Eq_OPM}
\Phi_{\bm{\tau}}(X,t)=\sum_{n\geq m_c+1} \Phi_n(\tau_n,\bflambda,X,t)  \boldsymbol{e}_n,
\ee
where either $ \Phi_n(\tau_n,\bflambda,X,t)$ is given by \eqref{Eq_param_timedependent} in the time-dependent case, or by \eqref{Eq_Phi_tau}, otherwise.  Whatever the case, the vector $\bm{\tau}^\ast$ is formed by the minimizers $\tau_n^\ast$ of $\mathcal{Q}_n$ given by \eqref{Eq_minQnHn}, for each $n$. Note that in the case of the time-dependent parameterization \eqref{Eq_param_timedependent}, the OPM reduced system  is run online by augmenting \eqref{Eq_OPM} with equations \eqref{Eq_for_I}, depending on the modes that are forced.

We emphasize finally that from a data-driven perspective, the OPM reduced system benefits from its dynamical origin. 
By construction, only a scalar parameter $\tau$ is indeed optimized per mode to parameterize. This parameter benefits furthermore from a dynamical interpretation since it plays a key role in balancing the small spectral gaps known as to be the main issue in applications of invariant or inertial manifold theory in practice \cite{zelik2014inertial}; see Remark \ref{Rem_smallgap} above.

\section{Predicting Tipping Points}\label{Sec_tippingmain}

\subsection{The Stommel-Cessi model and tipping points}
A simple model for oceanic circulation showing bistability is Stommel's box model \cite{stommel1961thermohaline}, where the ocean is represented by two boxes, a low-latitude box with temperature $T_1$ and salinity $S_1$, and a high-latitude box with temperature $T_2$ and salinity $T_2$; see  \cite[Fig.~1]{cessi1994simple}. The Stommel model can be viewed as the simplest ``thermodynamic'' version of the Atlantic Meridional Overturning Circulation (AMOC) \cite[Chap.~6]{dijkstra2013nonlinear}, a major ocean current system transporting warm surface waters toward the northern Atlantic that constitutes an important tipping point element of the climate system; see \cite{boers2021observation,lenton2008tipping}.

Cessi in  \cite{cessi1994simple} proposed a variation of this model, based on the box model of \cite{bryan1993toy}, consisting of an ODE system describing  the evolution of the differences $\Delta T = T_1-T_2$ and $\Delta S= S_1-S_2$; see \cite[Eq.~(2.3)]{cessi1994simple}.  The Cessi model trades the absolute functions involved in the original  Stommel model by polynomial relations  more prone to analysis.
The dynamics of  $\Delta T $ and  $\Delta S$ are coupled via the density difference  $\Delta \rho$, approximated by the relation $\Delta \rho=\alpha_S\Delta S-\alpha_T \Delta T$ which induces an exchange $Q$ of mass between the boxes to be given as $Q=1/\tau_d+(q/V) \Delta \rho^2$ according to Cessi's formulation, where $q$ denotes the the Poiseuille transport coefficient, $V$ the volume of a box, and $\tau_d$ the diffusion  timescale.   The coefficient $\alpha_S$ is a coefficient inversely proportional to the practical salinity unit, i.e.~unit based on the properties of sea water conductivity while  $\alpha_T$ is in $  ^{\circ}\,\textrm{C}^{-1}$; see \cite{cessi1994simple}. In this simple model, $\Delta T $ relaxes at a rate $\tau_r$ to a reference temperature $\theta$ (with $T_1=\theta/2$ and $T_2=-\theta/2$) in absence of coupling between the boxes.

Using the dimensionless variables $y=\alpha_S \Delta S/(\alpha_T \theta)$, $z=\Delta T/\theta$, and rescaling time by the diffusion timescale $\tau_d$,  the Cessi model can be written as \cite{berglund2002metastability}
\bea \label{Cessi_model}
& \dot{y} = F - y[1 + \mu(z - y)^2], \\
& \dot{z} = -\frac{1}{\epsilon}(z - 1) - z[1 + \mu(z - y)^2],
\eea
in which $F$ is proportional to the freshwater flux, $\epsilon=\tau_r/\tau_d$, and $\mu^2 = \tau_d (\alpha_T \theta)^2 q/V$; see \cite[Eq.~(2.6)]{cessi1994simple} and \cite{berglund2002metastability}.  

The nonlinear exchange of mass between the boxes is reflected by the nonlinear coupling terms in Eq.~\eqref{Cessi_model}. These nonlinear terms lead to multiple equilibria in certain parts of the parameter space, in particular when $F$ is varied over a certain range $[F_{c_1}, F_{c_2}]$, while $\mu$ and $\epsilon$ are fixed. One can even prove that Eq.~\eqref{Cessi_model} experiences two saddle-node bifurcations \cite{Kuznetsov04} leading to a typical S-shaped bifurcation diagram; see Fig.~\ref{Fig_Cessi_1st_result}A.

S-shaped bifurcation diagrams occur in oceanic models that go well beyond Eq.~\eqref{Cessi_model} such as based on the hydrostatic primitive equations or Boussinesq equations; see e.g.~\cite{thual1992catastrophe,dijkstra1997symmetry,weijer2003fully,weijer2019stability}.
More generally, S-shaped bifurcation diagrams and more complex multiplicity diagrams are known to occur in a broad class of nonlinear problems \cite{lions1982existence,kielhofer2011bifurcation,chekroun2018topological} that include energy balance climate models \cite{ghil1976climate, ghil1981energy,hetzer1997s,bodai2015global}, population dynamics models \cite{lee2011s,pruitt2018social},  vegetation pattern models \cite{bel2012gradual,zelnik2015gradual}, combustion models \cite{bebernes2013mathematical,frank2015diffusion,du2000exact,du2001proof}, and many other fields \cite{feudel2018multistability}.

The very presence of such S-shaped bifurcation diagrams provides the skeleton for tipping point phenomena  to take place when such models are subject to the appropriate stochastic disturbances and parameter drift, causing the system to ``tip'' or move away from one branch of attractors to another; see \cite{kuehn2011mathematical,ashwin2012tipping}. From an observational viewpoint, the study of tipping points has gained a considerable attention due to their role in climate change as a few components of the climate system (e.g.~Amazon forest, the AMOC) have been identified as candidates for experiencing such critical transitions if forced beyond the point of no return \cite{lenton2008tipping,caesar2018observed,boers2021observation}. 

Whatever the context, tipping  phenomena are due to a subtle interplay between nonlinearity, slow parameter drift, and fast disturbances. To better understand how these interactions cooperate to produce a tipping  phenomenon could help improve the design of early warning signals. Although, we will not focus on this latter point per se in this study, we show below, on the Cessi model, that the OPM framework provides useful insights in that perspective, by demonstrating its ability of deriving reduced models to predict accurately the crossing of a tipping point; see Sec.~\ref{Sec_tipping}.

\subsection{OPM results for a fixed $F$ value: Noise-induced transitions} 
\begin{figure*}[tbh!]
\centering
\includegraphics[width=\textwidth, height=0.35\textwidth]{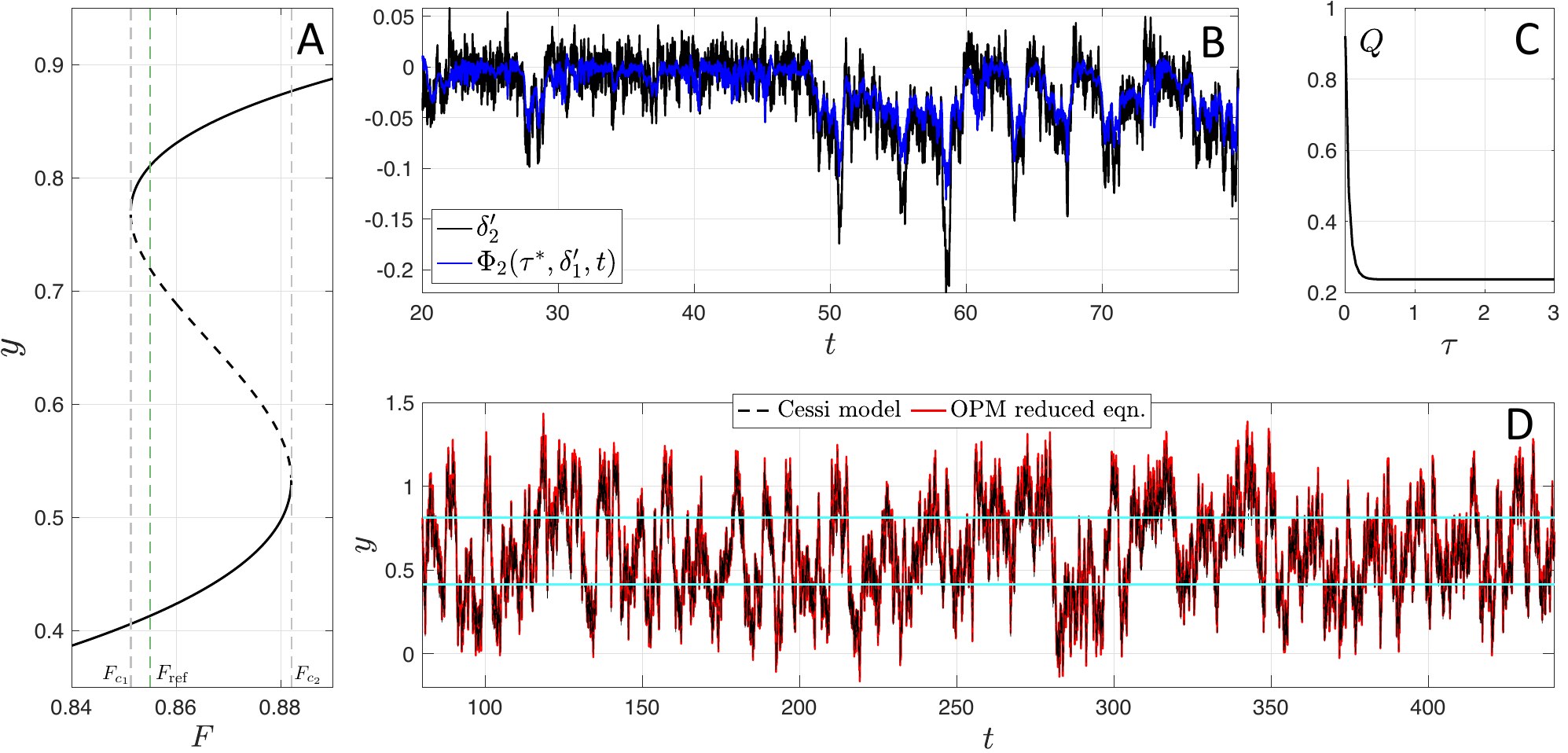}
\vspace{-2ex}
\caption{{\bf Panel A}: The S-shaped bifurcation of the Stommel-Cessi model~\eqref{Cessi_model} as the parameter $F$ is varied, shown here for $\mu = 6.2$ and $\epsilon = 0.1$. The two branches of locally stable steady states are marked by the solid black curves, and the other branch of unstable steady states are marked by the dashed black curve. The two vertical gray lines mark respectively $F_{c_1} = 0.8513$ and $F_{c_2} = 0.8821$, at which the saddle-node bifurcations occur. {\bf Panel B}: Parameterization of $\delta'_2$ by the OPM, $\Phi_2(\tau^*,\delta'_1,t)$ given by  \eqref{Eq_OPM_Cessi}, where $F$ is fixed to be $F_{\mathrm{ref}} = 0.855$ marked out by the vertical green line in Panel A and the noise strength parameter $\sigma$ in \eqref{Cessi_model_fluc2} is taken to be $\sigma = \sqrt{\epsilon}$, leading to $\sigma_1 \approx 0.3399$ and $\sigma_2 \approx -0.0893$ in \eqref{Eq_Cessi_eig}. {\bf Panel C}: The normalized parameterization defect $Q$ for $\Phi_2(\tau,\delta'_1,t)$ as $\tau$ is varied. {\bf Panel D}: The performance of the OPM reduced equation~\eqref{Eq_Cessi_reduced} in reproducing the noise-induced transitions experienced by $y(t)$ from the stochastically forced Stommel-Cessi model~\eqref{Cessi_model_fluc2}. Once \eqref{Eq_Cessi_reduced} is solved, the approximation of $y$ is constructed using \eqref{Eq_approx_y}. } \label{Fig_Cessi_1st_result}
\end{figure*}

We report in this section on the OPM framework skills to derive accurate reduced models to reproduce noise-induced transitions experienced by the Cessi model \eqref{Cessi_model}, when subject to fast disturbances, for a fixed value of $F$.
The training of the OPM operated here serves as a basis for the harder tipping point prediction problem dealt with in Sec.~\ref{Sec_tipping} below, 
when $F$ is allowed to drift slowly through the critical value $F_{c_2}$ at which the lower branch of steady states experiences a saddle-node bifurcation manifested by a turning point; see Fig.~\ref{Fig_Cessi_1st_result}A again. 
Recall that $F_{c_1}$ denotes here the $F$-value corresponding to the turning point experienced by the upper branch.

\noindent{\bf Reformulation of the Cessi model~\eqref{Cessi_model}.} 
The 1D OPM reduced equation is obtained as follows. First, we fix an arbitrary value of $F$ in $[F_{c_1}, F_{c_2}]$ that is denoted by $F_{\mathrm{ref}}$ and marked by the green vertical dashed line in Fig.~\ref{Fig_Cessi_1st_result}A.  
The  system \eqref{Cessi_model} is then rewritten for the fluctuation variables $y'(t) = y(t) - \overline{y}$ and $z'(t)= z(t) - \overline{z}$, where $\overbar{\bm{X}}=(\overline{y},\overline{z})$ denotes the steady state of Eq.~\eqref{Cessi_model} in the lower branch when $F=F_{\mathrm{ref}}$.

The resulting equation for $\bm{\delta}= (y', z')$ is then of the form 
\be \label{Cessi_model_fluc}
\dot{\bm{\delta}} = A \bm{\delta}+ G_2(\bm{\delta}) + G_3(\bm{\delta}),
\ee
with $A$, $G_2$, $G_3$ given by:
\begin{widetext}
\be\label{Eq_A}
A = \begin{pmatrix}
-1 - \mu (\overline{z} - \overline{y})^2 + 2 \mu \overline{y}(\overline{z} - \overline{y}) & -2 \mu \overline{y}(\overline{z} - \overline{y}) \\
2 \mu \overline{z}(\overline{z} - \overline{y}) & -\frac{1}{\epsilon} - 1 - \mu (\overline{z} - \overline{y})^2 - 2 \mu \overline{z}(\overline{z} - \overline{y})  
\end{pmatrix},
\ee
\be
G_2(\bm{\delta}) = \begin{pmatrix}
-\mu \overline{y}(z' - y')^2 - 2 \mu (\overline{z} - \overline{y}) y'(z' - y') \\
-\mu \overline{z}(z' - y')^2 - 2 \mu (\overline{z} - \overline{y}) z'(z' - y')    
\end{pmatrix} \quad \text{ and } \quad 
G_3(\bm{\delta}) = \begin{pmatrix}
-\mu y'(z' - y')^2 \\
-\mu z'(z' - y')^2
\end{pmatrix}. 
\ee
\end{widetext}

Since $\overbar{\bm{X}}$ is a locally stable steady state, we add noise to the first component of $\bm{\delta}$ to trigger transitions from the lower branch to the top branch in order to learn an OPM that can operate not only locally near $\overbar{\bm{X}}$ but also when the dynamics is visiting the upper branch. This leads to 
\be \label{Cessi_model_fluc2}
\dot{\bm{\delta}}= A \bm{\delta} + G_2(\bm{\delta}) + G_3(\bm{\delta}) + \bm{\sigma}\dot{W}(t),
\ee
where $\bm{\sigma}= (\sigma,0)^T$ and $W$ denotes a standard one-dimensional two-sided Brownian motion.
%

%%%%%%%%%%%%%%%%%%%%%%%%%%%%%%%%%%%%%%%%%%%%%%%%%%%%%%
%%%%%%%%%%%%%%%%%%%%%%%%%%%%%%%%%%%%%%%%%%%%%%%%%%%%%%

The eigenvalues of $A$ have negative real parts since $\overbar{\bm{X}}$ is locally stable.  
We assume that the matrix $A$ has two distinct real eigenvalues, which turns out to be the case for a broad range of parameter regimes. As in Sec.~\ref{Sec_LIAmain}, the spectral elements of the matrix $A$ (resp.~$A^{\ast}$) are denoted by $(\lambda_j, \boldsymbol{e}_j)$ (resp.~$(\lambda^*_j, \boldsymbol{e}^*_j)$), for $j=1,2$.  These eigenmodes are normalized to satisfy $\langle \boldsymbol{e}_j, \boldsymbol{e}^*_j  \rangle = 1$, and are bi-orthogonal otherwise.

Let us introduce $\bm{\delta}'=( \langle \bm{\delta}, \bm{e}^*_1 \rangle, \langle \bm{\delta}, \bm{e}^*_2 \rangle)^T $  and 
$\sigma_j=\langle\bm{\sigma}, \bm{e}^*_j \rangle$, for $j=1,2$. We also introduce  $\Lambda=\mbox{diag}(\lambda_1,\lambda_2)$. In the eigenbasis, Eq.~\eqref{Cessi_model_fluc2} can be written then as
\be \label{Eq_Cessi_eig}
 \dot{\bm{\delta}'}  = \Lambda \bm{\delta}'+ \mathcal{G}_2(\bm{\delta}') + \mathcal{G}_3(\bm{\delta}') + (\sigma_1\dot{W}(t),\sigma_2\dot{W}(t))^T,
\ee
with
\be \label{Eq_Gk_eigen_TJmodel}
\mathcal{G}_{k}^j(\bm{\delta}') = \Big \langle G_k \Big(\delta'_1 \bm{e}_{1} + \delta'_2 \bm{e}_{2} \Big), \bm{e}^*_j \Big \rangle,  \; j = 1,2,
\ee
where  $\delta'_j$ denotes the $j^{\mathrm{th}}$ component of $\bm{\delta}'$.

\noindent{\bf Derivation of the OPM reduced equation.}  We can now use the formulas of Sec.~\ref{Sec_explicit} to obtain an explicit variational parameterization  of the most stable direction carrying the variable $\delta'_2$ here, in terms of the least stable one, carrying $\delta'_1$.

For Eq.~\eqref{Eq_Cessi_eig}, both of the forcing terms appearing in the $p$- and $q$-equations of the corresponding BF system \eqref{Eq_BF_quad} are stochastic (and thus time-dependent). To simplify, we omit the stochastic forcing in Eq.~\eqref{BF1_quad} and work with the OPM formula  given by \eqref{Eq_param_timedependent} which, as shown below,  is sufficient for deriving an efficient reduced system.   

Thus, the formula \eqref{Eq_param_timedependent} becomes in this context
\bea \label{Eq_OPM_Cessi}
\Phi_2(\tau,X,t) & = D_{11}^2(\tau)B_{11}^2 X^2 + Z_{\tau}(t), 
\eea
where $D_{11}^2$ is given by \eqref{Eq_D_term0}, $B_{11}^2 = \langle \mathcal{G}_2(\bm{e}_1), \bm{e}^*_2 \rangle$, and 
\be \label{Eq_Z_term}
Z_{\tau}( t) = \sigma_2 e^{\lambda_2 t} \int_{t-\tau}^t e^{-\lambda_2 s} \dot{W}(s) \d s.
\ee

Once the optimal $\tau^*$ is obtained by minimizing the parameterization defect given by \eqref{Eq_minQnHn}, we obtain the following OPM reduced equation for $\delta'_1$:
\bea \label{Eq_Cessi_reduced}
\dot{X}  & = \lambda_1  X + \Big\langle G_{2}(X \bm{e}_1 + \Phi_2(\tau^*, X,t) \bm{e}_2), \bm{e}^*_1\Big\rangle   \\
& \quad +  \langle G_{3}(X\bm{e}_1 + \Phi_2(\tau^*, X,t) \bm{e}_2),  \bm{e}^*_1\rangle + \sigma_1 \dot{W}(t).
\eea
The online obtention of $X(t)$ by simulation of the OPM reduced equation Eq.~\eqref{Eq_Cessi_reduced} allows us to get the following approximation of the variables $(y(t),z(t))$ from the original Cessi model \eqref{Cessi_model}:
\be \label{Eq_approx_y}
(y_{\mathrm{app}}(t), z_{\mathrm{app}}(t))^T \hspace{-1ex}=  \hspace{-.5ex}X(t) \bm{e}_1  \hspace{-.25ex}+ \hspace{-.25ex}\Phi_2(\tau^*, X(t),t) \bm{e}_2 + \overbar{\bm{X}}, 
\ee
after going back to the original variables. 

\noindent{\bf Numerical results.} The OPM reduced equation Eq.~\eqref{Eq_Cessi_reduced} is able to reproduce the dynamics of $\delta'_1$ for a wide range of parameter regimes.  We show in Fig.~\ref{Fig_Cessi_1st_result} a typical example of skills for parameter values of $\mu$, $\epsilon$, $\sigma$ and $F=F_{\mathrm{ref}}$ as listed in Table \ref{Table_param_values}. 
Since $F_{\mathrm{ref}}$ lies in $[F_{c_1},F_{c_2}]$ (see Table \ref{Table_param_values}), the Cessi model \eqref{Cessi_model} admits three steady states  among which two are locally stable (lower and upper branches) and the other one is unstable (middle branch); see Fig.~\ref{Fig_Cessi_1st_result}A again.   For this choice of $F_{\mathrm{ref}}$, the steady state corresponding to the lower branch is $\overbar{X}=(\overline{y},\overline{z})$ with $\overline{y}=0.4130$, $\overline{z}=0.8285$, and the eigenvalues of $A$ are $\lambda_1 = -0.5168$, and $\lambda_2 = -15.7650$.

%%%%%%%%%%%%%%%%%%%%%
\begin{table}[ht] 
\caption{\small Parameter values} 
\label{Table_param_values}
\centering
\begin{tabular}{cccccccc}
\toprule\noalign{\smallskip}
& $\epsilon$ & $\sigma$ & $\mu$ & $F_{c_1}$ & $F_{c_2}$ & $F_{\mathrm{ref}}$   \\
\midrule\noalign{\smallskip}
& 0.1 & $\sqrt{\epsilon}$ & 6.2  & 0.8513 & 0.8821 & 0.855 \\
\noalign{\smallskip} \bottomrule 
\end{tabular}
\end{table}
%%%%%%%%%%%%%%%%%%%%%%%%%%

The offline trajectory $\bm{\delta}'(t)$ used as input for training $\Phi_2$ to find the optimal $\tau$, is taken as driven by an arbitrary Brownian path from Eq.~\eqref{Eq_Cessi_eig} for $t$ in the  time interval $[20, 80]$. 
The resulting offline skills of the OPM, $\Phi_2(\tau^\ast,\delta_1'(t),t)$, are shown as the blue curve in Fig.~\ref{Fig_Cessi_1st_result}B,  while the  original targeted time series $\delta'_2(t)$ to parameterize is shown in black. The optimal $\tau$ that minimizes the normalized parameterization defect $Q$ turns out to be $\infty$ for the considered regime, as shown in Fig.~\ref{Fig_Cessi_1st_result}C. One observes that the OPM captures, in average, the fluctuations of $\delta_2'(t)$; compare blue curve with black curve.  

The skills of the corresponding OPM reduced equation \eqref{Eq_Cessi_reduced} are shown in Fig.~\ref{Fig_Cessi_1st_result}D, after converting back to the original $(y,z)$-variable using \eqref{Eq_approx_y}. 
The results are shown out of sample, i.e.~for another noise path and over a time interval different from the one used for training. The 1D reduced  OPM reduced equation \eqref{Eq_Cessi_reduced} is able to capture the multiple noise-induced transitions occurring across the two equilibria (marked by the cyan lines), after transforming back to the original $(y,z)$-variable;  compare red with black curves in Fig.~\ref{Fig_Cessi_1st_result}D.
 Both the Cessi model \eqref{Cessi_model_fluc2} and the reduced OPM equation are numerically integrated using the Euler-Maruyama scheme with time step $\delta t = 10^{-3}$. 

Note that we chose here the numerical value of $F_{\mathrm{ref}}$ to be closer to $F_{c_1}$ than to $F_{c_2}$ for making  more challenging the tipping point prediction experiment conducted in Sec.~\ref{Sec_tipping} below. There, we indeed train the OPM for  $F=F_{\mathrm{ref}}$ while aiming at predicting the tipping phenomenon as  $F$ drifts slowly through $F_{c_2}$ (located thus further away) as time evolves.

%%%%%%%%%%%%%
\subsection{Predicting the crossing of tipping points}\label{Sec_tipping}

\noindent{\bf OPM reduced equation \eqref{Eq_Cessi_reduced} in the original coordinates.} For better interpretability, we first rewrite the OPM reduced equation~\eqref{Eq_Cessi_reduced} under the original coordinates in which the Cessi model \eqref{Cessi_model} is formulated. For that purpose, we exploit the components of the eigenvectors of $A$ given by \eqref{Eq_A} (for $F=F_{\mathrm{ref}}$) that we write as $\bm{e}_1 = (e_{11}, e_{12})^T$ and $\bm{e}_2 = (e_{21}, e_{22})^T$. Then, from \eqref{Eq_approx_y} the expression of  $y_{\mathrm{app}}$ rewrites as 
\bea \label{Eq_yapp}
y_{\mathrm{app}}(t) & = e_{11} X(t) + e_{21} \Phi_2(\tau^*, X(t),t) + \overline{y} \\
& = e_{11} X(t) + \gamma(\tau^*) X(t)^2 +  e_{21}Z_{\tau^*}(t) + \overline{y},
\eea
where the second line is obtained by using the expression of $\Phi_2$ given by \eqref{Eq_OPM_Cessi} and the notation $\gamma(\tau^*) = e_{21} D_{11}^2(\tau^*)B_{11}^2$. 
It turns out that 
\begin{widetext}
\be \label{Eq_X_using_y_app}
X(t) = \frac{- e_{11} + \sqrt{(e_{11})^2 - 4 \gamma(\tau^*) \big(e_{21} Z_{\tau^*}(t) + \overline{y} - y_{\mathrm{app}}(t)\big) }}{2 \gamma(\tau^*)} \stackrel{\mathrm{def}}{=} \varphi(\tau^*,y_{\mathrm{app}}(t), t), \; \; t \ge 0.
\ee  
\end{widetext}

From \eqref{Eq_approx_y}, we also have 
\bea \label{Eq_z_param}
z_{\mathrm{app}}(t) & = e_{12} X(t) + e_{22} \Phi_2(\tau^*, X(t),t) + \overline{z} \\
& = e_{12} X(t) + \gamma_2(\tau^*) X(t)^2 + e_{22} Z_{\tau^*}(t)  + \overline{z},
\eea
where $\gamma_2(\tau^*) = e_{22} D_{11}^2(\tau^*)B_{11}^2$. 

 We can now express $z_{\mathrm{app}}(t)$ as a function of $y_{\mathrm{app}}(t)$, i.e.~$z_{\mathrm{app}}(t) = \Psi(\tau^*,y_{\mathrm{app}}(t),t)$ with $\Psi$ given by 
 \bea \label{Eq_param_z}
\Psi(\tau^*,y_{\mathrm{app}},t) =\overline{z} &+ \gamma_2(\tau^*) \varphi(\tau^*,y_{\mathrm{app}}, t)^2  \\
& \;\; + e_{12} \varphi(\tau^*,y_{\mathrm{app}}, t) + e_{22} Z_{\tau^*}(t),
\eea
as obtained by replacing $X(t)$ with $\varphi(\tau^\ast,y_{\mathrm{app}}, t)$ in Eq.~\eqref{Eq_z_param}. 
Replacing $y_{\mathrm{app}}$ by the dummy variable $y$, one observes that  $\Psi$ provides a time-dependent manifold that parameterizes the variable  $z$ in the original Cessi model as (a time-dependent) polynomial function of radical functions in $y$  due to the expression of $\varphi(\tau^\ast,y, t)$; see \eqref{Eq_X_using_y_app}. 

We are now in position to derive the equation satisfied by $y_{\mathrm{app}}$ aimed at approximating $y$ in the original variables.  For that purpose, we first differentiate with respect to time the expression of $X(t)$ given in \eqref{Eq_X_using_y_app} to obtain 
\be \label{Eq_differential_relation}
\dot{X} = \frac{-e_{21} \dot{Z}_{\tau^*}(t) + \dot{y}_{\mathrm{app}}}{\sqrt{(e_{11})^2 - 4 \gamma(\tau^*) \big(e_{21} Z_{\tau^*}(t) + \overline{y} - y_{\mathrm{app}}\big)}}. 
\ee
On the other hand, by taking into account the expressions of $y_{\mathrm{app}}$ in \eqref{Eq_yapp} and  $z_{\mathrm{app}}$ in \eqref{Eq_z_param},  we observe that  the $X$-equation \eqref{Eq_Cessi_reduced} can be re-written  as 
\be \label{Eq_Cessi_reduced_v2}
\dot{X} = \langle \mathcal{F}(y_{\mathrm{app}}(t), z_{\mathrm{app}}(t))  + \bm{\sigma}\dot{W}(t), \bm{e}^*_1 \rangle.
\ee
with 
 \be
\mathcal{F}(y,z) = \begin{pmatrix}
&F_{\mathrm{ref}} - y[1 + \mu(z - y)^2] \\
&-\frac{1}{\epsilon}(z - 1) - z[1 + \mu(z - y)^2]
\end{pmatrix},
\ee
given by the right-hand side (RHS) of the Cessi model \eqref{Cessi_model} with $F = 
F_{\mathrm{ref}}$. 

Now by equating  the RHS of Eq.~\eqref{Eq_Cessi_reduced_v2} with that of Eq.~\eqref{Eq_differential_relation}, 
we obtain, after substitution of $z_{\mathrm{app}}$ by $\Psi(\tau^*,y_{\mathrm{app}},t)$, that $y_{\mathrm{app}}$ satisfies the following equation
\bea \label{Eq_Cessi_reduced3}
\dot{Y}   = \alpha(\tau^*,Y)  \Big\langle \mathcal{F}(Y,\Psi(\tau^*,Y,t)) &+ \bm{\sigma} \dot{W}(t), \bm{e}^*_1\Big\rangle  \\
& \qquad + e_{21} \dot{Z}_{\tau^*}(t),
\eea
where 
\be
\hspace{-.45ex}\alpha(\tau^*,Y) \hspace{-.2ex}= \hspace{-.2ex}\sqrt{(e_{11})^2 \hspace{-.1ex}- \hspace{-.05ex}4 \gamma(\tau^*) \big(e_{21} Z_{\tau^*}(t) + \overline{y} - Y\big)},
\ee
and $\dot{Z}_{\tau^*}$ is given by (cf.~\eqref{Eq_Z_term})
\be
\dot{Z}_{\tau^*}(t) = \lambda_2 Z_{\tau^*}(t) + \sigma_2 \dot{W}(t) - \sigma_2 e^{\lambda_2 \tau^*} \dot{W}(t -\tau^*).
\ee
Eq.~\eqref{Eq_Cessi_reduced3} is the OPM reduced equation for the $y$-variable of the Cessis model,  as rewritten 
in the original system of coordinates. This equation in its analytic formulation, mixes information about the two components of the RHS to Eq.~\eqref{Cessi_model} through the inner product involved in Eq.~\eqref{Eq_Cessi_reduced3}, while parameterizing the $z$-variable by the time-dependent parametarization $\Psi(\tau^*,\cdot,t)$ given by \eqref{Eq_param_z}. While designed for $F=F_{\mathrm{ref}}$ away from the lower tipping point, we show next that the OPM reduced system built up this way demonstrates a remarkable ability in predicting this tipping point when $F$ is allowed to drift away from $F_{\mathrm{ref}}$ in the course of time.

\noindent{\bf Prediction results.} Thus, we aim at predicting by our OPM reduced model, the tipping phenomenon experienced by the full model when $F$ is subject to a slow drift in time via
\be \label{Eq_F_drift}
F(t) =  F_0 + \kappa  t,
\ee
where $\kappa > 0$ is a small parameter, and $F_0$ is some fixed value of $F$ such that $F_0<F_{c_2}$, with $F_{c_2}$ denoting the parameter value of  the turning point of the lower branch of equilibria; see Fig.~\ref{Fig_Cessi_1st_result}A.

The original Cessi model is forced as follows   
\bea \label{Cessi_model_v3}
& \dot{y} = F(t) - y[1 + \mu(z - y)^2] + \sigma \dot{W}_t, \\
& \dot{z} = -\frac{1}{\epsilon}(z - 1) - z[1 + \mu(z - y)^2].
\eea

Introducing, $g(t)= \langle (F(t) - F_{\mathrm{ref}},0)^T, \bm{e}^*_1 \rangle$, and writing $F(t)$ as $F_{\mathrm{ref}} + (F(t) - F_{\mathrm{ref}})$, we consider the following forced version of the OPM reduced equation \eqref{Eq_Cessi_reduced3}: 
\bea \label{Eq_Cessi_reduced4}
\dot{Y}   = \alpha(\tau^*,Y)  \Big\langle \mathcal{F}(Y,\Psi(\tau^*,&Y,t)) + \bm{\sigma} \dot{W}(t), \bm{e}^*_1\Big\rangle  \\
& \;\; + e_{21} \dot{Z}_{\tau^*}(t) + g(t). 
\eea
 In this equation, recall that the parameterization, $\Psi(\tau^*,\cdot,t)$ has been trained for $F=F_{\mathrm{ref}}$; see above.

We set now $F_0 = 0.85$, $\kappa = 2\times 10^{-4}$, and $\sigma = \sqrt{\epsilon}/50$, while $\mu = 6.2$ and $\epsilon = 0.1$ are kept as in Table \ref{Table_param_values}.  Note that compared to the results shown in Fig.~\ref{Fig_Cessi_1st_result}, the noise level is here substantially reduced to focus on the tipping phenomenon to occur  in the vicinity of the turning point $F=F_{c_2}$.
The goal is thus to compare the behavior of $y(t)$ solving the forced Cessi model \eqref{Cessi_model_v3} with that of the solution $Y(t)$ of the (forced) OPM reduced equation \eqref{Eq_Cessi_reduced4}, as $F(t)$ crosses the critical value $F_{c_2}$.

\begin{figure}[tbh!]
\centering
\includegraphics[width=0.5\textwidth,height=.38\textwidth]{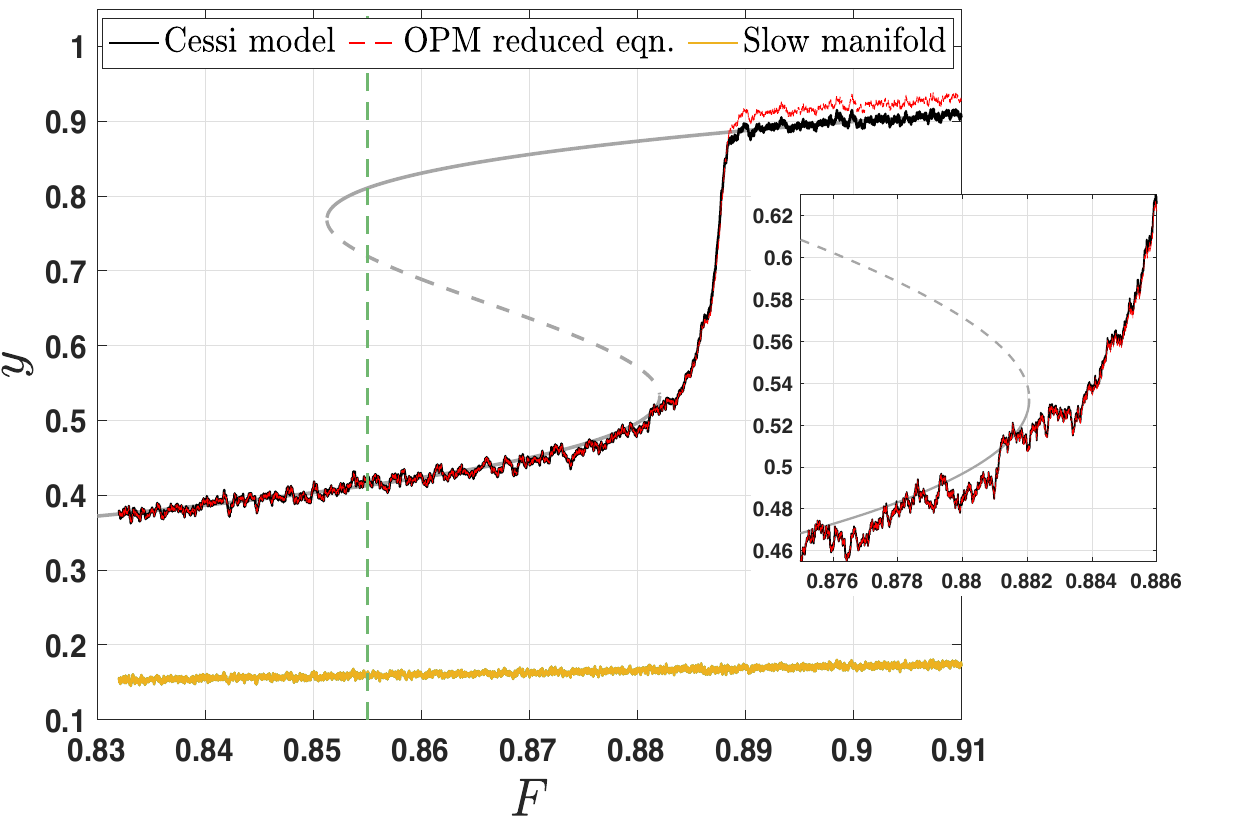}  
\vspace{-2ex}
\caption{The tipping phenomenon as predicted by the OPM reduced equation~\eqref{Eq_Cessi_reduced4} (red curve) compared with that experienced by the full Cessi model \eqref{Cessi_model_v3} (black curve). The OPM is trained for $F = F_{\mathrm{ref}} = 0.855$ as marked out by the vertical green line. Also shown in yellow is the result obtained from the slow manifold reduced Eq.~\eqref{Eq_slow}.} \label{Fig_Cessi_tipping_result}
\end{figure}

The results are shown in Fig.~\ref{Fig_Cessi_tipping_result}.  The red curve corresponds to the solution of the OPM reduced equation \eqref{Eq_Cessi_reduced4}, and the black curve to the $y$-component of the forced Cessi model \eqref{Cessi_model_v3}. 
Our baseline is the slow manifold parameterization which consists of simply parameterizing $z$ as $z=1+\mathcal{O}(\epsilon)$  in Eq.~\eqref{Eq_Cessi_reduced4} which provides, for $\epsilon$ sufficiently small due to the Tikhonov's theorem \cite{tikhonov1952systems,berglund2002metastability}, good reduction skills from the ``slow'' reduced equation, 
\be\label{Eq_slow}
\dot{y} = F(t) - y[1 + \mu(1 - y)^2] + \sigma \dot{W}_t.
\ee
Here, the value $\epsilon=0.1$ lies beyond the domain of applicability of the Tikhonov theorem, and as a result the slow reduced Eq.~\eqref{Eq_slow} fails in predicting any tipping  phenomenon; see yellow curve in Fig.~\ref{Fig_Cessi_tipping_result}.  In contrast, the OPM reduced equation \eqref{Eq_Cessi_reduced4} demonstrates a striking success in predicting the tipping phenomenon to the upper branch, in spite of being trained away from the targeted turning point, for $F=F_{\textrm{ref}}$ as marked by the (green) vertical dash line. The only caveat is the overshoot observed in the magnitude of the predicted random steady state in the upper branch.

%%%%%%%%%%%%%%%%%%%%%%%%%%%%%%%%%%%%%%%%%%%%%%%%%%%%%%
%%%%%%%%%%%%%%%%%%%%%%%%%%%%%%%%%%%%%%%%%%%%%%%%%%%%%%

The striking prediction result  of Fig.~\ref{Fig_Cessi_tipping_result} are shown for one noise realization. We explore now the accuracy in predicting such a tipping phenomenon  by the OPM reduced equation~\eqref{Eq_Cessi_reduced4}  when  the noise realization is varied. 

\begin{figure}[tbh!]
\centering
\includegraphics[width=0.48\textwidth,height=.32\textwidth]{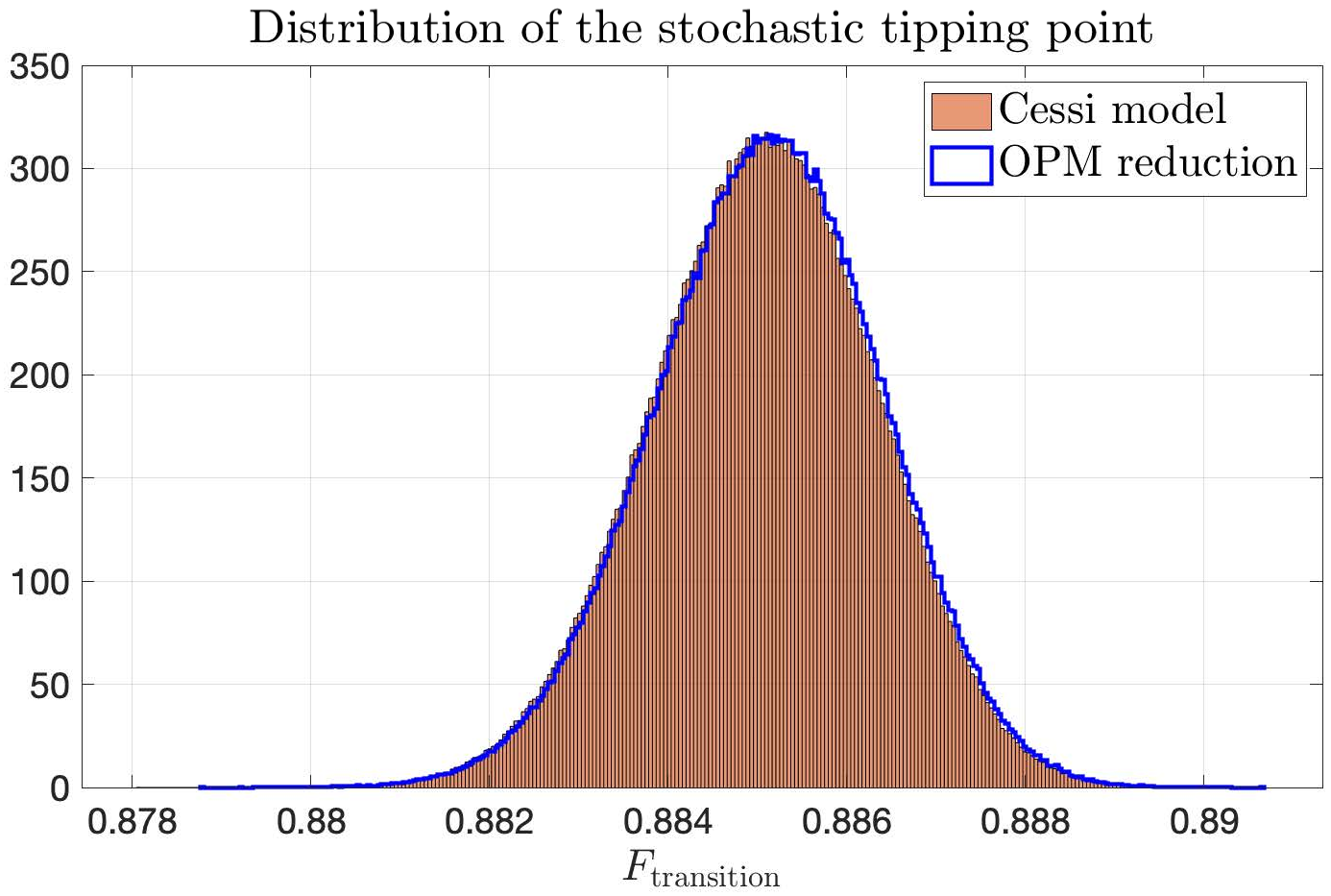} 
\vspace{-2ex}
\caption{Statistical distribution of the threshold value of $F$ at which the tipping phenomenon occurs for both the full Cessi model~\eqref{Cessi_model_v3} (orange bar) and the OPM reduced equation~\eqref{Eq_Cessi_reduced4} (blue curve). 
The histograms are computed based on $10^6$ arbitrarily fixed noise realizations.} \label{Fig_Cessi_stats}
\end{figure}

To do so, we estimate the statistical distribution of the $F$-value at which tipping takes place, denoted by $F_\mathrm{transition}$,  for  both the Cessi model \eqref{Cessi_model_v3} and Eq.~\eqref{Eq_Cessi_reduced4}. 
These distributions are estimated as follows. 
We denote by $\bar{y}_c$ the $y$-component of the steady state at  which the saddle-node bifurcation occurs in the lower-branch for $F=F_{c_2}$. As noise is turned on and $F(t)$ evolves slowly through $F_{c_2}$, the solution path $y(t)$ to Eq.~\eqref{Cessi_model_v3} increases in average while fluctuating around the  $\bar{y}_c$-value (due to small noise) before shooting off to the upper branch as  one nears $F_{c_2}$. During this process, there is a  time instant, denoted by $t_\mathrm{transition}$, such that for all $t > t_\mathrm{transition}$, $y(t)$  stay above $\bar{y}_c$. We denote the $F$-value corresponding to $t_\mathrm{transition}$ as $F_\mathrm{transition}$ according to \eqref{Eq_F_drift}. 
Whatever the noise realization, the solution to the OPM reduced equation~\eqref{Eq_Cessi_reduced4} experiences 
the same phenomenon leading thus to its own $F_\mathrm{transition}$ for a particular noise realization. 
As shown by the histograms in Fig.~\ref{Fig_Cessi_stats}, the distribution of $F_\mathrm{transition}$ predicted by the OPM reduced equation \eqref{Eq_Cessi_reduced4} (blue curve) follows closely that computed from the full system \eqref{Cessi_model_v3} (orange bars). 
Thus, not only the OPM reduced equation~\eqref{Eq_Cessi_reduced4} is qualitatively able to reproduce the passage through a tipping point (as shown in Fig.~\ref{Fig_Cessi_tipping_result}), but is also able to accurately predicting the critical $F$-value (or time-instant) at which the tipping phenomenon takes place with an overall success rate over 99\%.   
%%%%%%%%%%%%%%%%%%%%%%%%%%%%%%%%%%%%%%%%%%%%%%%%%%%%%%
%%%%%%%%%%%%%%%%%%%%%%%%%%%%%%%%%%%%%%%%%%%%%%%%%%%%%%

\section{Predicting Higher-Order Critical Transitions}\label{Sec_RB}

\subsection{Problem formulation}
In this section, we aim at applying our variational parameterization framework to the 
the following Rayleigh-B\'enard (RB) system from \cite{Reiterer_al98}:
%%%%%%%%%%%%%%%%%%%%%%%%%%%%%%%%%%%%
\bea \label{Eq_9DRBC}
\dot{C_1}&=-\sigma b_1 C_1 -C_2 C_4 +b_4 C_4^2+b_3 C_3 C_5 -\sigma b_2 C_7, \\
\dot{C_2}&=-\sigma C_2 +C_1 C_4 -C_2 C_5+ C_4 C_5 -\frac{\sigma}{2} C_9, \\
\dot{C_3}&=-\sigma b_1 C_3 +C_2 C_4 -b_4 C_2^2-b_3 C_1 C_5 +\sigma b_2 C_8, \\
\dot{C_4}&=-\sigma C_4 -C_2 C_3 - C_2C_5+C_4 C_5 +\frac{\sigma}{2}C_9, \\
\dot{C_5}&=-\sigma b_5 C_5 + \frac{1}{2} C_2^2 - \frac{1}{2} C_4^2, \\
\dot{C_6}&=- b_6 C_6 + C_2 C_9 - C_4C_9, \\
\dot{C_7}&=- b_1 C_7 -rC_1 + 2 C_5 C_8 - C_4 C_9, \\
\dot{C_8}&=- b_1 C_8  + rC_3 -2 C_5 C_7 +   C_2 C_9, \\
\dot{C_9}&=\hspace{-.1ex} - C_9\hspace{-.2ex}   -\hspace{-.2ex} C_2(r + 2  C_6 + C_8)\hspace{-.2ex}  +\hspace{-.2ex} C_4(r + 2  C_6+ C_7)
\eea
%%%%%%%%%%%%%%%%%%%%%%%%%%%%%%%%%%%%
Here $\sigma$ denotes the Prandtl number, and $r$ denotes the reduced Rayleigh number defined to be the ratio between the Rayleigh number $R$ and its critical value $R_c$ at which the convection sets in. The coefficients $b_i$ are given by
\beas
&b_1= \frac{4(1+a^2)}{1+2a^2},  \quad b_2=\frac{1+2a^2}{2(1+a^2)},\quad b_3=\frac{2(1-a^2)}{1+a^2}, \\
& b_4=\frac{a^2}{1+a^2}, \quad b_5=\frac{8a^2}{1+2a^2}, \quad  b_6=\frac{4}{1+2a^2},
\eeas
with $a=\frac{1}{2}$ corresponding to the critical horizontal wavenumber of the square convection cell. This system is obtained as a Fourier truncation of hydrodynamic equations describing Rayleigh-B\'enard (RB) convection in a 3D box \cite{Reiterer_al98}. 
 The Prandtl number is chosen to be $\sigma = 0.5$ in the experiments performed below, which is the same as used in \cite{Reiterer_al98}. The reduced Rayleigh number $r$ is varied according to these experiments; see Table \ref{Table_RB9D}.

Our goal is to assess the ability of our variational parameterization framework for predicting 
higher-order transitions arising in \eqref{Eq_9DRBC} by training the OPM only with data prior to the transition we aim at predicting.  
The Rayleigh number $r$ is our control parameter. As it increases, Eq.~\eqref{Eq_RBC_eigenbasis} undergoes several critical transitions/bifurcations, leading to chaos via a period-doubling cascade \cite{Reiterer_al98}. We focus on the prediction of two transition scenarios beyond Hopf bifurcation: (I) the period-doubling bifurcation, and (II) the transition from a period-doubling regime to chaos. Noteworthy are the failures  that standard invariant manifold theory or AIM encounter in the prediction of these transitions, here; see Appendix \ref{Sec_IMfailue}.  

To do so, we re-write Eq.~\eqref{Eq_9DRBC} into the following compact form: 
\be
\dot{\boldsymbol{C}} = L \boldsymbol{C} + B(\boldsymbol{C},\boldsymbol{C}),
\ee
where $\boldsymbol{C}=(C_1, \cdots C_9)^T$, and $L$ is the $9\times9$ matrix given by
\bes
L= \hspace{-.5ex}\begin{pmatrix}
-\sigma b_1 & 0 & 0 & 0 & 0 & 0 & -\sigma b_2 & 0 & 0 \\  % row 1
0 & -\sigma & 0 & 0 & 0 & 0 & 0 & 0 & -\frac{\sigma}{2}  \\ % row 2
0 & 0 &-\sigma b_1& 0 & 0 & 0 & 0 & \sigma b_2 &  0  \\ % row 3
0 & 0 & 0 & -\sigma & 0 & 0  & 0 & 0 & \frac{\sigma}{2} \\ % row 4
0 & 0 & 0 & 0 & -\sigma b_5 & 0 &  0  & 0 & 0 \\ % row 5
0 & 0 & 0 & 0 & 0 & - b_6 & 0  & 0  & 0 \\ % row 6
-r & 0 & 0 & 0 & 0 & 0 & - b_1  & 0  & 0 \\ % row 7
0 & 0 & r & 0 & 0 & 0 & 0 & - b_1  & 0 \\ % row 8
0 & -r & 0 & r & 0 & 0 & 0 & 0  & -1 % row 9
\end{pmatrix}
\ees
The nonlinear term $B$ is defined as follows.  For any $\boldsymbol{\phi}=(\phi_1, \cdots, \phi_9)^{T}$ and $\boldsymbol{\psi}=(\psi_1, \cdots, \psi_9)^{T}$ in $\mathbb{C}^9$, we have
\bea\label{B_term_RBC}
B(\boldsymbol{\phi},\boldsymbol{\psi}) =  \begin{pmatrix}
-\phi_2 \psi_4 +b_4 \phi_4 \psi_4 +b_3 \phi_3 \psi_5 \\
\phi_1 \psi_4 -\phi_2 \psi_5+ \phi_4 \psi_5 \\
\phi_2 \psi_4 -b_4 \phi_2\psi_2-b_3 \phi_1 \psi_5 \\
-\phi_2 \psi_3 - \phi_2\psi_5+\phi_4 \psi_5 \\
\frac{1}{2} \phi_2 \psi_2 - \frac{1}{2} \phi_4 \psi_4 \\
\phi_2 \psi_9 - \phi_4\psi_9 \\
2 \phi_5 \psi_8 - \phi_4 \psi_9 \\
-2 \phi_5 \psi_7 +   \phi_2 \psi_9 \\
-2 \phi_2 \psi_6 - \phi_2 \psi_8 +2 \phi_4 \psi_6+  \phi_4 \psi_7 
\end{pmatrix}.
\eea
We now re-write   Eq.~\eqref{Eq_9DRBC} in terms of fluctuations with respect to its mean state. 
In that respect, we subtract from $\boldsymbol{C}(t)= (C_1(t), \cdots, C_9(t))$ its mean state $\overline{\boldsymbol{C}}_r$, which is estimated, in practice, from simulation of Eq.~\eqref{Eq_9DRBC} over a typical characteristic time of the dynamics that resolves e.g.~decay of correlations (when the dynamics is chaotic) or a period (when the dynamics is periodic); see also \cite{CLM19_closure}. 
The corresponding ODE system for the fluctuation variable, $\boldsymbol{\delta}(t) = \boldsymbol{C}(t) - \overline{\boldsymbol{C}}_r$, is then given by:
\be \label{Eq_RBC_fluct}
\frac{\mathrm{d} \boldsymbol{\delta}}{\mathrm{d} t} = A \boldsymbol{\delta} + B(\boldsymbol{\delta},\boldsymbol{\delta}) + L \overline{\boldsymbol{C}}_r + B(\overline{\boldsymbol{C}}_r, \overline{\boldsymbol{C}}_r),
\ee
with 
\be\label{RB_Linear_part_perturbed}
A \boldsymbol{\delta}=L \boldsymbol{\delta} + B(\overline{\boldsymbol{C}}_r,\boldsymbol{\delta}) + B(\boldsymbol{\delta},\overline{\boldsymbol{C}}_r). 
\ee
Denote the spectral elements of the matrix $A$ by $\{(\lambda_j, \boldsymbol{e}_j) \; : \; 1 \le j \le 9\}$ and those of $A^{\ast}$ by $\{(\lambda^*_j, \boldsymbol{e}^*_j) \; : \; 1 \le j \le 9\}$.  Here the eigenmodes are normalized so that $\langle \boldsymbol{e}_j, \boldsymbol{e}^*_k  \rangle = 1$ if $j=k$, and $0$ otherwise. By taking the expansion of $\boldsymbol{\delta}$ in terms of the eigenelements of $L$, we get
\be
\boldsymbol{\delta}(t)= \sum_{j=1}^9 y_j (t) \boldsymbol{e}_j \quad \text{ with } \quad y_j (t) = \langle\boldsymbol{\delta}(t), \boldsymbol{e}^*_j \rangle.
\ee
Assuming that $A$ is diagonalizable in $\mathbb{C}$ and by rewriting \eqref{Eq_RBC_fluct} in the variable $\boldsymbol{y}=(y_1, \cdots, y_9)^{T}$, 
 we obtain that 
\bea \label{Eq_RBC_eigenbasis}
&\dot{y}_j = \lambda_j(\overline{\boldsymbol{C}}_r) y_j + \sum_{k,\ell = 1}^9  B_{k \ell}^j(\overline{\boldsymbol{C}}_r)  y_k y_{\ell} +  F_j(\overline{\boldsymbol{C}}_r), \\
& j = 1, \cdots, 9,
\eea
where $B_{k \ell}^j(\overline{\boldsymbol{C}}_r)=\langle B(\boldsymbol{e}_k, \boldsymbol{e}_\ell), \boldsymbol{e}^*_j\rangle$ and 
\be\label{Eq_forcing}
F_j(\overline{\boldsymbol{C}}_r) = \langle L \overline{\boldsymbol{C}}_r + B(\overline{\boldsymbol{C}}_r, \overline{\boldsymbol{C}}_r), \boldsymbol{e}^*_j\rangle.
\ee  
Like the $\lambda_j$, the eigenmodes also depend on the mean state $\overline{\boldsymbol{C}}_r$, explaining the dependence of the interaction coefficients  $B_{k \ell}^j$. From now on, we work with Eq.~\eqref{Eq_RBC_eigenbasis}.

%%%%%%%%%%%%%%%%%%%%%%%%%
\begin{figure*}
\centering
\includegraphics[width=1\linewidth, height=.5\textwidth]{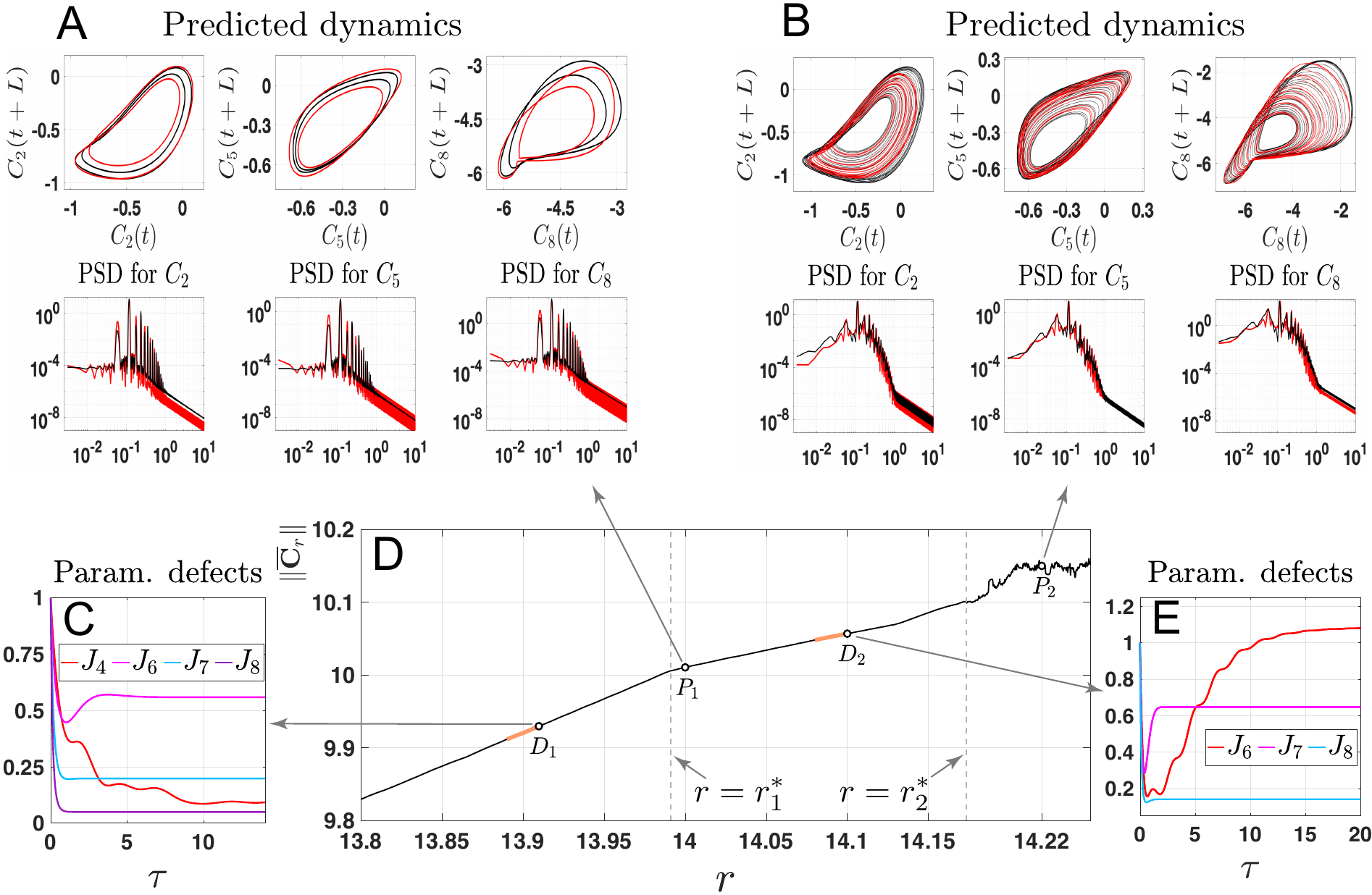}
\caption{{\bf Prediction of transitions: OPM prediction resulls}. {\bf Panel D:} The black curve shows the dependence on $r$ of the norm  of the mean state vector, $\overline{\boldsymbol{C}}_r$. The vertical dashed lines at $r_1^* = 13.991$ and $r_2^* = 14.173$ mark the onset of period-doubling bifurcation and chaos, respectively. The points $P_1$ and $P_2$ correspond to the $r$-values at which the two prediction experiments are conducted; see Table \ref{Table_RB9D}. The orange segments that precede the points $D_1$ and $D_2$, denote the parameter intervals over which data is used to build the OPM reduced system for predicting the dynamics at $r=r_{P_1}$ and $r=r_{P_2}$, respectively; see Steps 1-4.  
The  (normalized) parameterization defects are shown in {\bf Panels C} and {\bf E} for training data at $r=r_{D_1}$ and $r=r_{D_2}$, respectively. 
 {\bf Panels A and B:} Global attractor (in lagged coordinates) and PSDs for three  selected components at $r=r_{P_1}$ and $r=r_{P_2}$, respectively. Black curves are from the full system \eqref{Eq_9DRBC} and red ones from the OPM reduced system \eqref{Eq_RBC_reduced}.}
\label{Fig_RB_9D_combo}
\end{figure*}

 \subsection{Predicting Higher-order Transitions via OPM: Method}\label{Sec_predict_method}

Thus, we aim at predicting for Eq.~\eqref{Eq_RBC_eigenbasis} two types of transition: (I) the period-doubling bifurcation, and (II) the transition from period-doubling to chaos, referred to as Experiments I and II, respectively.
For each of these experiments, we are given a reduced state space $H_\c=\mbox{span}(\boldsymbol{e}_1,\cdots, \boldsymbol{e}_{m_c})$ with $m_c$ as indicated in Table \ref{Table_RB9D}, depending on the parameter regime.   The eignemodes are here  ranked by the real part of the corresponding eigenvalues, $\bm{e}_1$ corresponding to the eigenvalue with the largest real part.  
The goal is to derive an OPM reduced system \eqref{Eq_reduced} able to predict such transitions. The challenge lies in the optimization stage of the OPM as due to the prediction constraint, one is prevented to use data from the full model at the parameter $r=r_P$ which one desires to predict the dynamics.  Only data prior to  the critical value $r^*$ at which the concerned transition,  either period-doubling or chaos, takes place, are here allowed.

%%%%%%%%%%%%%%%%%%%%%%%%%%%%%%%%%%%%%%%%%%%%%%%%
\begin{table}
\centering 
\caption{{\bf Prediction experiments for the RB System \eqref{Eq_RBC_eigenbasis}}.
The parameter values $r_D (<r^\ast)$  corresponds to the allowable upper bound for which the mean state dependence of $\overline{\boldsymbol{C}}_r$ on $r$ is estimated, in view of extrapolation at $r=r_P$. In each experiment, $I_r=[r_0,r_D]$ with $r_D-r_0=2\times 10^{-2}$, corresponding to segments show in orange in Fig.~\ref{Fig_RB_9D_combo}-(D).   The critical value $r^\ast$ indicates the parameter value at which the transition occurs, depending on the experiment. The parameter value  $r_P>r^\ast$ at which the prediction is sought, is also given.}
\begin{tabular}{ccccc}
 & $m_c$ & $r_D$ & $r^\ast$ & $r_P$ \\
\midrule
Experiment I (period-doubling)& 3 &   $13.91$ & $13.99$ &  $14$  \\
Experiment II  (chaos)& 5 & $14.10$ & $14.17$ &  $14.22$ \\
\bottomrule
\end{tabular} \label{Table_RB9D}
%%%%%%%%%%%%%%
\end{table}
%%%%%%%%%%%%%%%%%%%%%%%%%%%%%%%%%%%%%%%%%%%%%%%%

 Due to the dependence on $\overline{\boldsymbol{C}}_r$ of $\boldsymbol{\lambda}$, as well as of the interaction coefficients $B_{k \ell}^j$ and forcing terms $F_j$, a particular attention to this dependence has to be paid. Indeed, recall that the parameterizations $\Phi_n$ given by the explicit formula \eqref{Eq_Phi_tau},
 depend heavily here on the spectral elements of linearized operator $A$ at $\overline{\boldsymbol{C}}_r$, and thus does  its optimization. Since the goal is to eventually dispose of an OPM reduced system able to predict the dynamical behavior at $r=r_P>r^\ast$, one cannot rely on data from the full model at  $r=r_P$ and thus one cannot exploit in particular the knowledge of  the mean state $\overline{\boldsymbol{C}}_r$ for $r=r_P$. We are thus forced to estimate the latter for our purpose.

To do so, we estimate the dependence on $r$ of $\overline{\boldsymbol{C}}_r$ on an interval $I_r= [r_0, r_D]$ such that $r_0<r_D<r^\ast$ (see Table \ref{Table_RB9D}), and use this estimation to extrapolate the value of $\overline{\boldsymbol{C}}_r$ at $r=r_P$, that we denote by $\overline{\boldsymbol{C}}^{\text{ext}}_{r_P}$.  For both experiments, it turned out that a linear extrapolation is sufficient. 
This extrapolated value $\overline{\boldsymbol{C}}^{\text{ext}}_{r_P}$ allows us to compute the spectral elements of  the operator $A$ given by \eqref{RB_Linear_part_perturbed} in which  $\overline{\boldsymbol{C}}^{\text{ext}}_{r_P}$ replaces the genuine mean state $\overline{\boldsymbol{C}}_r$. Obviously, the better is the approximation  of $\overline{\boldsymbol{C}}_r$ by   $\overline{\boldsymbol{C}}^{\text{ext}}_{r_P}$, the better the approximation of $ \boldsymbol{\lambda}(\overline{\boldsymbol{C}}_{r_P})$ by $\boldsymbol{\lambda}(\overline{\boldsymbol{C}} ^{\text{ext}}_{r_P})$ along with the corresponding eigenspace.

We turn now to the training of the OPM exploiting these surrogate spectral elements, that will be used for predicting the dynamical behavior at $r=r_P$. To avoid any looking ahead whatsoever, we allow ourselves to only use the data of the full model \eqref{Eq_RBC_eigenbasis} for $r=r_D<r_P$ and minimize the following 
parameterization defect 
\be\label{Eq_cost_func2}
\mathcal{J}_n(\tau)=\overline{\big|y_n^{r_D}(t) - \Phi_n(\tau,\boldsymbol{\lambda}(\overline{\boldsymbol{C}}^{\text{ext}}_{r_P}), \boldsymbol{y}_\c^{r_D} (t)) \big|^2},
\ee  
for each relevant $n$, in which $\boldsymbol{y}_\c^{r_D} (t)$ (resp.~$y_n (t)$) denotes the full model solution's  amplitude in $H_\c$ (resp.~$\boldsymbol{e}_n$) at $r=r_D$. The normalized defect is then given by $J_n(\tau)=\mathcal{J}_n(\tau)/\overline{|y_n^{r_D}|^2}$. To compute these defects we proceed as follows. We first use the spectral elements  at $r=r_P$ of the linearized matrix $A$ at the extrapolated mean state $\overline{\boldsymbol{C}}^{\text{ext}}_{r_P}$ to form the coefficients $R_n$, $\lambda_n^{-1}(1 - e^{\tau \lambda_n}) F_n$ and $D_{ij}^n(\tau, \boldsymbol{\lambda}) B_{ij}^n$  involved in the expression  \eqref{Eq_Phi_tau} of $\Phi_n$. For each $n$, the defect $\mathcal{J}_n(\tau)$ is then fed by the input data at $r=r_D$ (i.e.~$y_n^{r_D}(t)$ and $\boldsymbol{y}_\c^{r_D}(t)$), and evaluated for $\tau$ that  varies in some sufficiently large interval to capture global minima. The results are shown in Fig.~\ref{Fig_RB_9D_combo}C and  Fig.~\ref{Fig_RB_9D_combo}E for Experiment I and II, respectively. Note that the evaluation of $\mathcal{J}_n(\tau)$ for a range of $\tau$ is used here only for visualization purpose  as it is not needed to find a minimum. A simple gradient-descent algorithm can be indeed used for the latter; see \cite[Appendix]{CLM19_closure}.  We denote by $\wh{\tau}_n^\ast$ the resulting optimal value of $\tau$ per variable to parameterize.

After the  $\wh{\tau}_n^\ast$ are found for each $n\geq m_c+1$,
the OPM reduced system used to predict the dynamics at $r=r_P$ takes the form:
\begin{widetext}
\bea \label{Eq_RBC_reduced}
\dot{X}_j = \wh{\lambda}_j(r_P) X_j & + \sum_{k,\ell = 1}^{m_c} \wh{B}_{k \ell}^j(r_P) X_k X_\ell 
  +  \sum_{k = 1}^{m_c} \sum_{\ell=m_c+1}^9 \big(\wh{B}_{k \ell}^j(r_P) + \wh{B}_{\ell k}^j(r_P) \big) X_k \Phi_\ell(\wh{\tau}_\ell^\ast,\wh{\boldsymbol{\lambda}}(r_P),X) \\
&  +  \sum_{k,\ell=m_c+1}^9  \wh{B}_{k \ell}^j(r_P) \Phi_k(\wh{\tau}_k^\ast, \wh{\boldsymbol{\lambda}}(r_P), X)  \Phi_\ell(\wh{\tau}_\ell^\ast, \wh{\boldsymbol{\lambda}}(r_P),X)  + \wh{F}_j(r_P), \quad j = 1, \cdots, m_c,
\eea
with $\wh{\lambda}_j(r_P)=\lambda_j(\overline{\boldsymbol{C}}^{\text{ext}}_{r_P})$, $\wh{B}_{k \ell}^j(r_P)=B_{k \ell}^j(\overline{\boldsymbol{C}}^{\text{ext}}_{r_P})$, and $\wh{F}_j(r_P)=F_j(\overline{\boldsymbol{C}}^{\text{ext}}_{r_P}).$
\end{widetext}
%%%%%%%%%%%%%%%%

We can then summarize  our approach for predicting higher-order transitions via OPM reduced systems as follows:
\bi
\item[{\bf Step 1.}] Extrapolation $\overline{\boldsymbol{C}}^{\text{ext}}_{r_P}$ of $\overline{\boldsymbol{C}}_r$ at $r=r_P$ (the parameter at which one desires to predict the dynamics).
\item[{\bf Step 2.}] Computation of the spectral elements of the linearized operator $A$ at $r=r_P$, by replacing  $\overline{\boldsymbol{C}}_r$ by $\overline{\boldsymbol{C}}^{\text{ext}}_{r_P}$ in \eqref{RB_Linear_part_perturbed}.
\item[{\bf Step 3.}] Training of the OPM using the spectral elements of Step 2 and data of the full model for $r=r_D<r_P$.  
\item[{\bf Step 4.}] Run the OPM reduced system \eqref{Eq_RBC_reduced} to predict the dynamics at $r=r_P$.

\ei

We mention that the minimization of certain parameterization defects may require a special care such as for $J_6$ in Experiment II.
Due to the presence of nearby local minima (see red curve in Fig.~\ref{Fig_RB_9D_combo}E), the analysis of the optimal value of $\tau$ to select for calibrating an optimal parameterization of $y_6(t)$ is more subtle and exploits actually a complementary metric known as the parameterization correlation; see Sec.~\ref{Sec_localmin}.

Obviously, the accuracy in approximating the genuine mean state $\overline{\boldsymbol{C}}_{r_P}$  by $\wh{\boldsymbol{C}}_{r_P}$ is a determining factor in the transition prediction procedure described in Steps 1-4 above. Here, the relative error in approximating $\|\overline{\boldsymbol{C}}_{r_{P_i}}\|$ ($i=1,2$) is $0.03\%$ for Experiment I, and  $0.27\%$ for Experiment II. For the latter case, although the parameter dependence is rough beyond $r_2^\ast$ (see Panel D), there is no brutal local variations of relative large magnitude as identified for other nonlinear systems \cite{Chek_al14_RP,chekroun2017emergence}. Systems for which a linear response to parameter variation  is a valid assumption \cite{gallavotti1995dynamical,lucarini2018revising}, thus constitute seemingly a favorable ground to apply the proposed  transition prediction procedure.

%%%%%%%%%%%%%%%%%%%%%%%%%%%%%%%%%%%%%%%%%%%%%%%%
\subsection{Prediction of Higher-order Transitions by OPM reduced systems}\label{Sec_OPMpredictionRB}

 As summarized in Figs.~\ref{Fig_RB_9D_combo}A and B, for each prediction experiment of Table~\ref{Table_RB9D}, the OPM reduced system \eqref{Eq_RBC_reduced} not only successfully predicts the occurrence of the targeted transition at either $r=r_{P_1}$ and $r=r_{P_2}$ ($P_1$ and $P_2$ in  Fig.~\ref{Fig_RB_9D_combo}D), but also approximates with good accuracy the embedded (global) attractor as well as key statistics of the dynamics such as the power spectral density (PSD). Note that for Experiment I (period-doubling) we choose the reduced dimension to be $m_c=3$ and to be $m_c=5$ for Experiment II (transition to chaos). For Experiment I the global minimizers of
$\mathcal{J}_n(\tau)$ given by  \eqref{Eq_cost_func2} are retained to build up the corresponding OPM. In Experiment II, all the global minimizer  are  also retained except for $\mathcal{J}_6(\tau)$ from which 
 the second minimizer is used for the construction of the OPM.

Once the OPM is built, an approximation of $\boldsymbol{C}(t)$ is obtained from the solution $X(t)$ to the  corresponding OPM reduced system \eqref{Eq_RBC_reduced}  according to 
\bea \label{Eq_C_approx}
\boldsymbol{C}^{\text{PM}}(t)  = \overline{\boldsymbol{C}}^{\text{ext}}_{r_P}  +& \sum_{j=1}^{m_c} X_j(t) \boldsymbol{e}_j \\
 &+ \sum_{n=m_c+1}^{N}\Phi_n(\wh{\tau}_n^\ast, \wh{\boldsymbol{\lambda}}(r_P),X) \boldsymbol{e}_n,
\eea
with $\Phi_n$ given by \eqref{Eq_Phi_tau} whose spectral entries are given by the spectral elements of $A$ 
given by \eqref{RB_Linear_part_perturbed} for $r=r_P$ in which $\overline{\bm{C}}_{r_P}$ is replaced by $ \overline{\boldsymbol{C}}^{\text{ext}}_{r_P}$. The $\boldsymbol{e}_j$ and $\bm{e}_n$ denote the eigenvectors of this matrix and $N=9$ here.

As baseline, the skills of the OPM reduced system \eqref{Eq_RBC_reduced} are compared to those from reduced systems when parameterizations from invariant manifold theory such as \eqref{Eq_LP_integral}  or from  inertial manifold theory such as \eqref{Eq_FMTgen}, are employed; see \cite[Theorem 2]{CLM19_closure} and Remark \ref{Rem_FMT}. The details and results are discussed in Appendix \ref{Sec_IMfailue}. The main message is that compared to the OPM reduced system \eqref{Eq_RBC_reduced}, the reduced systems based on these traditional inertial/inertial manifold theories fail in predicting the correct dynamics. A main player in this failure lies in the inability of these standard parameterizations to accurately approximate small-energy variables that are dynamically important; see Appendix \ref{Sec_IMfailue}. The OPM parameterization by its variational nature enables to fix this over-parameterization issue here.

\subsection{Distinguishing between close local minima: Parameterization correlation analysis}\label{Sec_localmin} 
Since the $\Phi_n$'s coefficients depend nonlinearly on $\tau$ (see Eqns.~(\ref{Eq_Phi_tau})-(\ref{Eq_D_term0}) and \eqref{Formula_GammaF}), the parameterization defects, $\mathcal{J}_n(\tau)$, defined in \eqref{Eq_cost_func2}  are also highly nonlinear, and may not be convex in the $\tau$-variable as shown for certain variables in Panels C and E of Fig.~\ref{Fig_RB_9D_combo}. A minimization algorithm  to reach most often its global minimum is nevertheless detailed in \cite[Appendix A]{CLM19_closure}, and is not limited to the RB system analyzed here. In certain rare occasions, a local minimum may be an acceptable halting point with 
an online performance slightly improved compared to that of the global minimum. 
In such a case, one discriminates between a local minimum and the global one by typically inspecting another metric offline: the correlation angle that measures essentially the collinearity between the actual high-mode dynamics and its parameterization. Here, such a situation occurs for Experiment II; see $J_6$ in Fig.~\ref{Fig_RB_9D_combo}E.

Following \cite[Sec.~3.1.2]{CLM19_closure}, we  recall thus a simple criterion to assist the selection of an OPM when there are more than one local minimum displayed by the parameterization defect and the corresponding local minimal values are close to each other.  

Given a parameterization $\Phi$ that is not trivial (i.e.~$\Phi\neq 0$), we define the {\it parameterization correlation} as, 
\be\label{Eq_corr_param}
c(t)= \frac{\mathrm{Re}\langle \Phi(y_{\c}(t)), y_{\s}(t) \rangle}{\|\Phi(y_{\c}(t))\| \; \|y_{\s}(t)\|}.
\ee
It provides a measure of collinearity between the unresolved variable $y_{\s}(t)$ and its parameterization $\Phi(y_{\c}(t))$, as time evolves.  It serves thus as a complimentary, more geometric way to measure the phase coherence between the resolved and unresolved variables than with the parameterization defect $\mathcal{Q}_n(\tau_n, T)$ defined in \eqref{Eq_minQnHn}. The closer to unity $c(t)$ is for a given parameterization, the better the phase coherence between the resolved and unresolved variables is expected to hold.

\begin{figure}[h]
\centering
\includegraphics[width=0.45\textwidth, height=0.25\textwidth]{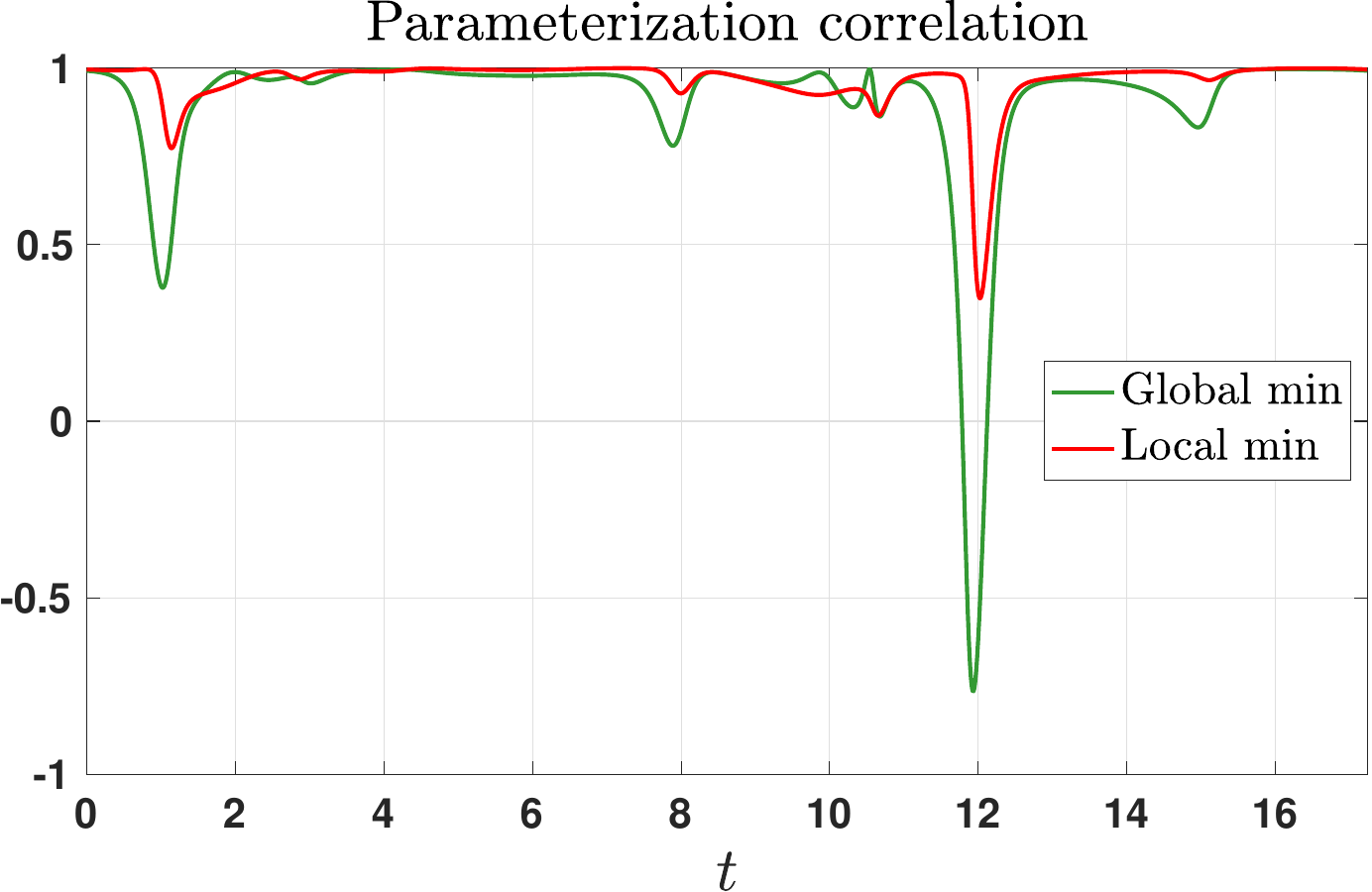}
\vspace{-2ex}
\caption{{\bf Parameterization correlation $c(t)$ defined by \eqref{Eq_corr_param} for the OPM in Experiment II (transition to chaos)}. The metric $c(t)$ is here computed by using the full model solutions for $r=r_D=14.1$ in $J_6$ in \eqref{Eq_cost_func2}, by making $\tau=0.65$ and $\tau=1.80$ corresponding to the global minimizer (green curve) and the nearby local minimizer (red curve), respectively; see Fig.~\ref{Fig_RB_9D_combo}E.} \label{Fig_RB9D_Global_vs_local}
\end{figure}

We illustrate this point on $J_6(\tau)$ shown in Fig.~\ref{Fig_RB_9D_combo}E. 
The defect $J_6$ exhibits two close minima corresponding to $J_6\approx 0.1535$ and $J_6\approx 0.1720$, occurring respectively at $\tau \approx 0.65$ and $\tau\approx1.80$. Thus, the parameterization defect alone does not help provide a satisfactory discriminatory diagnosis between these two minimizers. To the contrary,  the parameterization correlation shown in Fig.~\ref{Fig_RB9D_Global_vs_local} allows for diagnosing more clearly that the OPM associated with the local minimizer has a neat advantage compared  to that associated with the global minimizer.

%%%%%%%%%%%%%%%%%%%%%%%%%%%%%%%%%%%%%%%%%%%%%%%%%%%%%%
%%%%%%%%%%%%%%%%%%%%%%%%%%%%%%%%%%%%%%%%%%%%%%%%%%%%%%

As explained in \cite[Sec.~3.1.2]{CLM19_closure}, this parameterization correlation criterion is also useful to assist the selection of the dimension of the reduced state space. Indeed, once an OPM has been determined, the dimension $m_c$ of the reduced system should be chosen so that the parameterization correlation is sufficiently close to unity as measured for instance in the $L^2$-norm over the training time window. The basic idea is that one should not only parameterize properly the statistical effects of the neglected variables but also avoid to lose their phase relationships with the unresolved variables \cite{mccomb2001conditional,CLM19_closure}. For instance, for predicting the transition to chaos, we observe that an OPM comes with a parameterization correlation much closer to unity in the case  $m_c=5$ than $m_c=3$ (not shown).

%%%%%%%%%%%%%%%%%%%%%%%%%%%%%%%%%%%%%%%%%%%%%%%%%%%%%%%
%%%                                              
%%%         
%%%%%%%%%%%%%%%%%%%%%%%%%%%%%%%%%%%%%%%%%%%%%%%%%%%%%%%

\section{Concluding Remarks} \label{Sec_conclusion}

In this article, we have described a general framework, based on the
the OPM approach introduced in \citep{CLM19_closure} for autonomous systems, to derive effective  (OPM) reduced systems able to either predict  higher-order transitions caused by parameter regime shift, or  tipping phenomenas caused by model's stochastic  disturbances and slow parameter drift.  In each case, the OPMs are sought as continuous deformations of classical invariant manifolds to handle parameter regimes away from bifurcation onsets while keeping the reduced state space relatively low-dimensional. 
The underlying OPM parameterizations are derived analytically from model's equations, constructed as explicit solutions to auxiliary BF systems such as Eq.~\eqref{Eq_BF_withforcing}, whose backward integration time is optimized per unresolved mode. This optimization involves the minimization of  natural cost functionals tailored upon the full dynamics at  parameter values prior the transitions take place; see \eqref{Eq_minQnHn}.  In each case---either for prediction of higher-order transitions or tipping points---we presented compelling evidences of success of the OPM approach to address such extrapolation problems typically difficult to handle for data-driven reduction methods.

As reviewed in Sec.~\ref{Sec_LIAmain}, the BF systems such as Eq.~\eqref{Eq_BF_withforcing} allow for drawing  insightful relationships  with the approximation theory of classical invariant/inertial manifolds (see \cite[Sec.~2]{CLM19_closure} and \cite{haro2016parameterization,PR06,PR09}) when the backward integration time $\tau$ in Eq.~\eqref{Eq_BF_withforcing} is sent to $\infty$; see Theorem~\ref{Thm_BF}. 
Once the invariant/inertial manifolds fail to provide legitimate parameterizations, the optimization of the backward integration time may lead to local minima at finite $\tau$ of the parameterization defect which if sufficiently small, give access to skillful  parameterizations to predict e.g.~higher-order transitions. This way, as illustrated in Sec.~\ref{Sec_RB} (Experiment II), the OPM approach allows us to bypass the well-known stumbling blocks tied to the presence of small spectral gaps such as encountered in inertial manifold theory \cite{zelik2014inertial}, and tied to the accurate parameterization of dynamically important small-energy variables; see also Remark~\ref{Rem_smallgap}, Appendix \ref{Sec_IMfailue} and \cite[Sec.~6]{CLM19_closure}.

To understand the dynamical mechanisms at play behind a tipping phenomenon and to predict its occurrence are of uttermost importance, but this task is hindered by the often high dimensionality nature of the underlying physical system.
Devising accurate, and analytic reduced models able to predict tipping phenomenon from complex system is thus of prime importance to serve understanding.
The OPM results of Sec.~\ref{Sec_tippingmain} indicate great promises for  the OPM approach to tackle this important task for more complex systems. Other reduction methods to handle the prediction of tipping point dynamics have been proposed recently in the literature but mainly for mutualistic networks \cite{jiang2018predicting,jiang2019harnessing,laurence2019spectral}. The OPM  formulas presented here are not limited to this kind of networks and can be directly applied to a broad class of spatially extended systems governed by (stochastic) PDEs or to nonlinear time-delay systems by adopting the framework of \cite{CGLW16} for the latter to devise center manifold parameterizations and their OPM generalization; see  \cite{chekroun2020efficient,Chekroun_al22SciAdv}.

The OPM approach could be also informative to design for such systems nearing a tipping point, early warning signal (EWS) indicators from multivariate time series. Extension of EWS techniques to multivariate data is an active field of research with methods ranging from empirical orthogonal functions reduction \cite{held2004detection} to methods exploiting the system's Jacobian matrix and relationships with the cross-correlation matrix \cite{williamson2015detection} or exploiting the detection of spatial location of  ``hotspots'' of stability loss \cite{bathiany2013detecting,kefi2014early}. By its nonlinear modus operandi for reduction, the OPM parameterizations identifies subgroup of essential variables to characterize a tipping phenomenon which in turn could be very useful to construct the relevant observables of the system for the design of useful EWS indicators.

Thus, since the examples of Secns.~\ref{Sec_tippingmain} and \ref{Sec_RB} are representative of more general  problems of prime importance, the successes demonstrated by the OPM approach on these examples invite for further investigations for more complex systems.

Similarly, we would like to point out another important aspect in that perspective.  For spatially extended systems, the modes involve typically wavenumbers that can help interpret certain primary physical patterns. Formulas such as \eqref{Eq_Phi_tau} or \eqref{Eq_param_timedependent} arising in OPM parameterizations involve a rebalancing of the interactions $B_{ij}^n$ among such modes by means of the coefficients $D_{ij}^n(\tau, \boldsymbol{\lambda})$; see Remark~\ref{Rem_smallgap}. Thus,  an OPM parameterization when skillful to derive efficient systems may provide useful insights to explain emergent patterns due to certain nonlinear interactions of wavenumbers, for regimes beyond the instability onset.  For more complex systems than dealt with here, it is known already than near the instability onset of primary bifurcations, center manifold theory provides such physical insights; see e.g.~\cite{chekroun2022transitions,dijkstra2015dynamic,han2020dynamic,hsia2008attractor,hsia2010rotating,lu2019hopf,sengul2015pattern}.

By the approach proposed here relying on OPMs obtained as continuous deformations of invariant/center manifolds, one can thus track the interactions that survive or emerge between certain wavenumbers as one progresses through higher-order transitions. Such insights  bring new elements to potentially explain the low-frequency variability of recurrent large-scale patterns typically observed e.g.~in oceanic models \cite{berloff1999large,dijkstra2005low},  offering at least new lines of thoughts to the dynamical system approach proposed in previous studies \cite{nadiga2001global,simonnet2003low,Simonnet2005,dijkstra2005low}. In that perspective, the analytic expression \eqref{Formula_GammaF} of terms such as $R_n(F, \bm{\lambda},\tau , X)$ in \eqref{Eq_Phi_tau} bring also new elements to break down the nonlinear interactions between the forcing of certain modes compared to others. 
Formulas such as \eqref{Formula_GammaF} extend to the case of  time-dependent or stochastic forcing of  the coarse-scale modes whose exact expression will be communicated elsewhere. As for the case of noise-induced transitions reported in Fig.~\ref{Fig_Cessi_1st_result}D, these generalized formulas are expected to provide in particular new insights to regime shifts not involving crossing a bifurcation point  but tied to other mechanisms such as slow-fast cyclic transitions  or stochastic resonances \cite{dakos2015resilience}. 

Yet, another important practical aspect that deserve in-depth investigation is related to the robustness of the OPM approach subject to model errors. Such errors can arise from multiple sources, including e.g.~imperfection of the originally utilized high-fidelity model in capturing the true physics and noise contamination of the solution data used to train the OPM. For the examples of Secns.~\ref{Sec_tippingmain} and \ref{Sec_RB}, since we train the OPM in a parameter regime prior to the occurrence of the concerned transitions, the utilized training data contain, de facto, model errors. The reported results in these sections show that the OPM is robust in providing effective reduced models subject to such model errors. Nevertheless, it would provide useful insights if one can systematically quantify uncertainties of the OPM reduced models subject to various types of model errors. In that respect,  several approaches could be beneficial for stochastic systems such as those based on the information-theoretic framework of \cite{majda2018model} or the perturbation theory of ergodic Markov chains and linear response theory \cite{zhang2021error}, as well as methods based on data-assimilation techniques \cite{gottwald2021supervised}.

\begin{acknowledgements}
We are grateful to the reviewers for their constructive comments, which helped enrich the discussion and presentation. This work has been partially supported by  the Office of Naval Research (ONR) Multidisciplinary University Research Initiative (MURI) grant N00014-20-1-2023, by the National Science Foundation grant DMS-2108856, and by  the European Research Council (ERC) under the European Union's Horizon 2020 research and innovation program (grant agreement no. 810370).  
We also acknowledge the computational resources provided by Advanced Research Computing at Virginia Tech for the results shown in Fig.~\ref{Fig_Cessi_stats}. 
\end{acknowledgements}

\section*{Data Availability}
The data that support the findings of this study are available
from the corresponding author upon reasonable request.

\appendix
{\small 
\section{Proof of Lemma \ref{Lemma_foundation}}\label{Proof_Lemma}

We introduce the notations
\bea
&g_\tau(t)=\int_{-\tau}^0 e^{-s A_{\s}} \Pi_\s F(s+t)\d s,\\
& \Psi_{\tau}(X)=\int_{-\tau}^0 e^{-s A_{\s}} \Pi_{\s} B(e^{s A_{\c}} X) \d s. 
\eea

In Eq.~\eqref{eq_h1_NDS}, replacing $X$ by $e^{r A_{\c}}X$ where $r$ is an arbitrary real number, we get
\be
\Phi_\tau (e^{r A_{\c}}X,t) = \Psi_{\tau}( e^{r A_{\c}} X)  + g_\tau(t). 
\ee
Note that
\bea
\varphi(r) &= \Psi_{\tau}( e^{r A_{\c}} X) = \int_{-\tau}^0 e^{-s A_{\s}} \Pi_{\s} B(e^{ (s+r) A_{\c}} X) \d s \\
& = \int_{r-\tau}^{r} e^{(r-s') A_{\s}} \Pi_{\s} B(e^{s' A_{\c}} X) \d s'. 
\eea
We obtain then
\bea \label{dPhi_eq1}
\partial_r \varphi(r) = \Pi_{\s} & B( e^{r A_{\c}} X)\\
& - e^{\tau A_{\s}}  \Pi_{\s} B( e^{(r-\tau) A_{\c}} X) + A_{\s} \varphi(r).
\eea
On the other hand, since $\varphi(r) =\Phi_\tau (e^{r A_{\c}} X,t) - g_\tau(t)$ due to \eqref{eq_h1_NDS}, we also have
\be \label{dPhi_eq2}
\partial_r \varphi(r)   = D \Phi_\tau (e^{r A_{\c}}X,t) A_{\c} e^{r A_{\c}} X.
\ee
By taking the limit $r \rightarrow 0$, we obtain from \eqref{dPhi_eq1}--\eqref{dPhi_eq2} that 
\bea \label{Eq_invariance_NDS_partI}
D \Phi_\tau (X,t) A_{\c} X = \Pi_{\s} B(X)  - e^{\tau A_{\s}}  & \Pi_{\s} B( e^{-\tau A_{\c}} X) \\
& + A_{\s} \Psi_{\tau}(X).
\eea

Now if we replace $t$ in \eqref{eq_h1_NDS} by $t+r$, we obtain
\be
\Phi_\tau (X,t+r) = \Psi_{\tau}(X)  + g_\tau(t+r). 
\ee
Note that
\bea
\widetilde{\varphi}(r) & = g_\tau(t+r)  =  \int_{-\tau}^0 e^{-s A_{\s}} \Pi_\s F(s+t+r)\d s \\
& = \int_{r-\tau}^r e^{(r-s') A_{\s}} \Pi_\s F(s'+t)\d s'. 
\eea
It follows that
\be \label{dPhi_eq3}
\partial_r \widetilde{\varphi}(r) = \Pi_{\s} F(r+t)  - e^{\tau A_{\s}}  \Pi_{\s} F( r-\tau+t) + A_{\s} \widetilde{\varphi}(r).
\ee
Since $\widetilde{\varphi}(r) =\Phi_\tau (X,t+r) - \Psi_{\tau}(X)$, we also have
\be \label{dPhi_eq4}
\partial_r \widetilde{\varphi}(r)  = \partial_t \Phi_\tau (X,t+r). 
\ee
By taking the limit $r \rightarrow 0$, we obtain from \eqref{dPhi_eq3}--\eqref{dPhi_eq4} that 
\bea \label{Eq_invariance_NDS_partII}
\partial_t \Phi_\tau (X,t) =  \Pi_{\s} &F(t)  \\
&- e^{\tau A_{\s}}  \Pi_{\s} F(t-\tau) + A_{\s} g_\tau(t).
\eea
The homological equation \eqref{Eq_invariance_NDS} follows then  from  \eqref{Eq_invariance_NDS_partI} and \eqref{Eq_invariance_NDS_partII} while gathering the relevant terms to form the operator $\mathcal{L}_A$ given by \eqref{Def_L}.

\section{The low-order term in the OPM parameterization}\label{Sec_LIA_formula}
%%%%%%%%%%%%%%%%%%%%%%
In the special case where $F$ in \eqref{Eq_BF_quad} is time-independent, by introducing $X_i = \langle X, \boldsymbol{e}^*_i \rangle$ and $F_i = \langle F, \boldsymbol{e}^*_i  \rangle$,  the $R_n$-term in the parameterization
\eqref{Eq_Phi_tau}  takes the form
\bea \label{Formula_GammaF}
R_n(F,\bflambda,\tau,X)  &=  \sum_{i, j = 1}^{m_c} U_{ij}^n(\tau, \boldsymbol{\lambda}) B_{ij}^n F_{i}F_{j}\\
& + \sum_{i, j = 1}^{m_c}   V_{ij}^n(\tau, \boldsymbol{\lambda})  F_{j} (B^n_{ij} + B^n_{ji})   X_{i},
\eea
\begin{widetext}
with 
\be\label{Eq_Un}
U_{ij}^n(\tau, \boldsymbol{\lambda}) = \begin{cases}
 \frac{1}{\lambda_{i} \lambda_{j}}\Big(D_{ij}^n(\tau, \boldsymbol{\lambda})  
 - \frac{1 - \exp(-\tau(\lambda_{i} - \lambda_{n}))}{\lambda_{i} - \lambda_{n}} \\
\qquad\qquad - \frac{1 - \exp(-\tau(\lambda_{j} - \lambda_{n}))}{\lambda_{j} - \lambda_{n}} 
 - \frac{1 - \exp(\tau \lambda_{n})}{\lambda_{n}} \Big), & \text{if $\lambda_{i} \neq 0$ and $\lambda_{j} \neq 0$}, \vspace{2ex}\\
 \frac{1}{\lambda_{i}}\Big(\frac{\tau \exp(-\tau(\lambda_{i} - \lambda_{n}))}{\lambda_{i} - \lambda_{n}} - \frac{1 - \exp(-\tau(\lambda_{i} - \lambda_{n}))}{(\lambda_{i} - \lambda_{n})^2} \\
\qquad\qquad + \frac{\tau \exp(\tau \lambda_{n})}{\lambda_n} + \frac{1-\exp(\tau \lambda_{n})}{(\lambda_{n})^2} \Big), & \text{if $\lambda_{i} \neq 0$ and $\lambda_{j} = 0$}, \vspace{2ex}\\
 \frac{1}{\lambda_{j}}\Big(\frac{\tau \exp(-\tau(\lambda_{j} - \lambda_{n}))}{\lambda_{j} - \lambda_{n}} - \frac{1 - \exp(-\tau(\lambda_{j} - \lambda_{n}))}{(\lambda_{j} - \lambda_{n})^2} \\
\qquad\qquad + \frac{\tau \exp(\tau \lambda_{n})}{\lambda_n} + \frac{1-\exp(\tau \lambda_{n})}{(\lambda_{n})^2} \Big), & \text{if $\lambda_{i} = 0$ and $\lambda_{j} \neq 0$}, \vspace{2ex}\\
\frac{\tau^2 \exp(\tau \lambda_{n})}{\lambda_{n}} - \frac{2}{\lambda_n} \Big(
\frac{\tau \exp(\tau \lambda_{n})}{ \lambda_{n}} + \frac{1-\exp(\tau \lambda_{n})}{(\lambda_{n})^2} \Big), & \text{if $\lambda_{i} = 0$ and $\lambda_{j} = 0$},
 \end{cases}
\ee
and 
\be\label{Eq_Vn}
V_{ij}^n(\tau, \boldsymbol{\lambda}) = \begin{cases}
 \frac{1 - \exp(-\tau(\lambda_{i} + \lambda_{j} - \lambda_{n}))}{\lambda_{j}(\lambda_{i} + \lambda_{j} - \lambda_{n})} - \frac{1 - \exp(-\tau(\lambda_{i} - \lambda_{n}))}{\lambda_{j}(\lambda_{i} - \lambda_{n})}, & \text{if $\lambda_{j} \neq 0$}, \vspace{2ex}\\
 \frac{\tau \exp(-\tau(\lambda_{i} - \lambda_{n}))}{\lambda_{i} - \lambda_{n}} - \frac{1 - \exp(-\tau(\lambda_{i} - \lambda_{n}))}{(\lambda_{i} - \lambda_{n})^2}. & \text{otherwise}.
 \end{cases}
\ee
\end{widetext}
%%%%%%%%%%%%%%%%%%

For the application to the Rayleigh-B\'enard system considered in Sec.~\ref{Sec_RB}, the forcing term $F$ is produced after rewriting the original unforced system in terms of the fluctuation variable with respect to the mean state. For this problem, the eigenvalues $\boldsymbol{\lambda}$ and the interaction coefficients $B_{ij}^n$ both depend on the mean state.
%%%%%%%%%%%%%%%%%%%%%%%%%%%%%%%%%%%%%%%%%%%%%%%%%%%%%%

\section{Numerical aspects}
In the numerical experiments of Sec.~\ref{Sec_RB},  the full RB system \eqref{Eq_9DRBC} as well as the OPM reduced system~\eqref{Eq_RBC_reduced} are numerically integrated using a standard fourth-order Runge--Kutta (RK4) method with a time-step size taken to be $\delta t = 5 \times 10^{-3}$. Note that since the eigenvalues of $A$ are typically complex-valued, some care is needed when integrating \eqref{Eq_RBC_reduced}. Indeed, since complex eigenmodes of $A$ must appear in complex conjugate pairs, the corresponding components of $\boldsymbol{x}$ in \eqref{Eq_RBC_reduced} form also complex conjugate pairs. To prevent round-off errors that may disrupt the underlying complex conjugacy, after each RK4 time step, we enforce complex conjugacy as follows. Assuming that $x_j $ and $x_{j+1}$ form a conjugate pairs and that after an RK4 step, $x_j = a_1 + i b_1$ and $x_{j+1} = a_2 - i b_2$  ($i^2=-1$) where $a_1, a_2, b_1$ and $b_2$ are real-valued with $a_1\approx a_2$ and $b_1\approx a_2$, we redefine $x_{j}$ and $x_{j+1}$ to be respectively given by $x_j = (a_1 + a_2)/2 + i (b_1+b_2)/2$ and $x_{j+1} = (a_1 + a_2)/2 - i (b_1+b_2)/2$. For each component corresponding to a real eigenmode, after each RK4 time step, we simply redefine it to be its real part and ignore the small imaginary part that may also arise from round-off errors. The same numerical procedure  is adopted to integrate the reduced systems obtained from invariant manifold reduction as well as the FMT parameterization used in Appendix \ref{Sec_IMfailue} below.

\section{Failure from standard manifold parameterizations}\label{Sec_IMfailue}
In this section we briefly report on the reduction skills for the chaotic case of Sec.~\ref{Sec_OPMpredictionRB} as   achieved through application of the invariant manifold theory or standard formulas used in inertial manifold theory. 
Invariant manifold theory is usually applied near a steady state. To work within a favorable ground for these theories, we first chose the steady state $\overline{\boldsymbol{Y}}_{\! r}$ that is closest to the mean state $\overline{\boldsymbol{C}}_r$ at the parameter value $r$ we want to approximate the chaotic dynamics via a reduced system. This way the chaotic dynamics is located near this steady state.

To derive the invariant manifold (IM) reduced system, we also re-write the RB system \eqref{Eq_9DRBC} in  the fluctuation variable, but this time using the steady state $\overline{\boldsymbol{Y}}_{\! r}$, namely for the $\boldsymbol{\delta}(t) = \boldsymbol{C}(t) - \overline{\boldsymbol{Y}}_{\! r}$.
The analogue of Eq.~\eqref{Eq_RBC_fluct} is then projected onto the eigenmodes of linearized part 
\be\label{Eq_A_forFMT}
\tilde{A} \boldsymbol{\delta}=L \boldsymbol{\delta} + B(\overline{\boldsymbol{Y}}_{\! r}, \boldsymbol{\delta}) + B(\boldsymbol{\delta},\overline{\boldsymbol{Y}}_{\! r}),
\ee
to obtain: 
\bea \label{Eq_RBC_eigenbasis2}
\dot{w}_j = \lambda_j w_j + \sum_{k,\ell = 1}^9  \tilde{B}_{k \ell}^j  w_k w_{\ell}, \;\; j = 1, \cdots, 9.
\eea
Here the $(\lambda_j, \boldsymbol{f}_j)$ denote the spectral elements of $\tilde{A}$, while the interaction coefficients are given as $\tilde{B}_{k \ell}^j=\langle B(\boldsymbol{f}_k, \boldsymbol{f}_\ell), \boldsymbol{f}^*_j\rangle$ with $\boldsymbol{f}^*_j$ denoting the eigenvector of $\tilde{A}^*$ associated with $\lambda_j^*$. Note that unlike \eqref{Eq_RBC_eigenbasis}, there is no constant forcing term on the RHS of \eqref{Eq_RBC_eigenbasis2} since $L \overline{\boldsymbol{Y}}_{\! r} + B(\overline{\boldsymbol{Y}}_{\! r}, \overline{\boldsymbol{Y}}_{\! r}) = 0$ as $\overline{\boldsymbol{Y}}_{\! r}$ is a steady state.  

Here also,  the reduced state space  $H_\c=\mbox{span}(\boldsymbol{e}_1,\cdots, \boldsymbol{e}_{m_c})$ with $m_c=5$ as indicated in Table \ref{Table_RB9D}, for the chaotic regime. The  local IM associated with $H_{\c}$ has its components approximated by \cite[Theorem 2]{CLM19_closure}: 
\be \label{Eq_h2}
h_n(\boldsymbol{X}) = \sum_{i, j = 1}^{m_c} \frac{\tilde{B}_{ij}^n}{\lambda_{i} + \lambda_{j} - \lambda_{n}} X_{i} X_{j}, \quad n = m_c+1, \cdots, 9.
\ee
when $\Re(\lambda_{i} + \lambda_{j} - \lambda_{n}) >0$ for $1\le i,j\le m_c$ and $n \ge m_c+1$.

The corresponding IM reduced system takes the form given by \eqref{Eq_RBC_reduced} in which the components, $\Phi_\ell$, of the OPM therein are replaced by the $h_\ell$ given by \eqref{Eq_h2} for $\ell=m_c+1, \cdots, 9$.  Unlike the OPM reduced system \eqref{Eq_RBC_reduced}, here we use the true eigenvalues $\lambda_j$'s  at the parameter $r$ at which the dynamics is chaotic ($r=r_{P_2}$ in the notations of Fig.~\ref{Fig_RB_9D_combo}).

%%%%%%%%%%%%%%%%%%%%%%%%%%%%%%%%%%%%%%%%%%%%%%%%%%%%%%
\begin{figure}[tbh!]
\centering
\includegraphics[width=0.48\textwidth, height=0.3\textwidth]{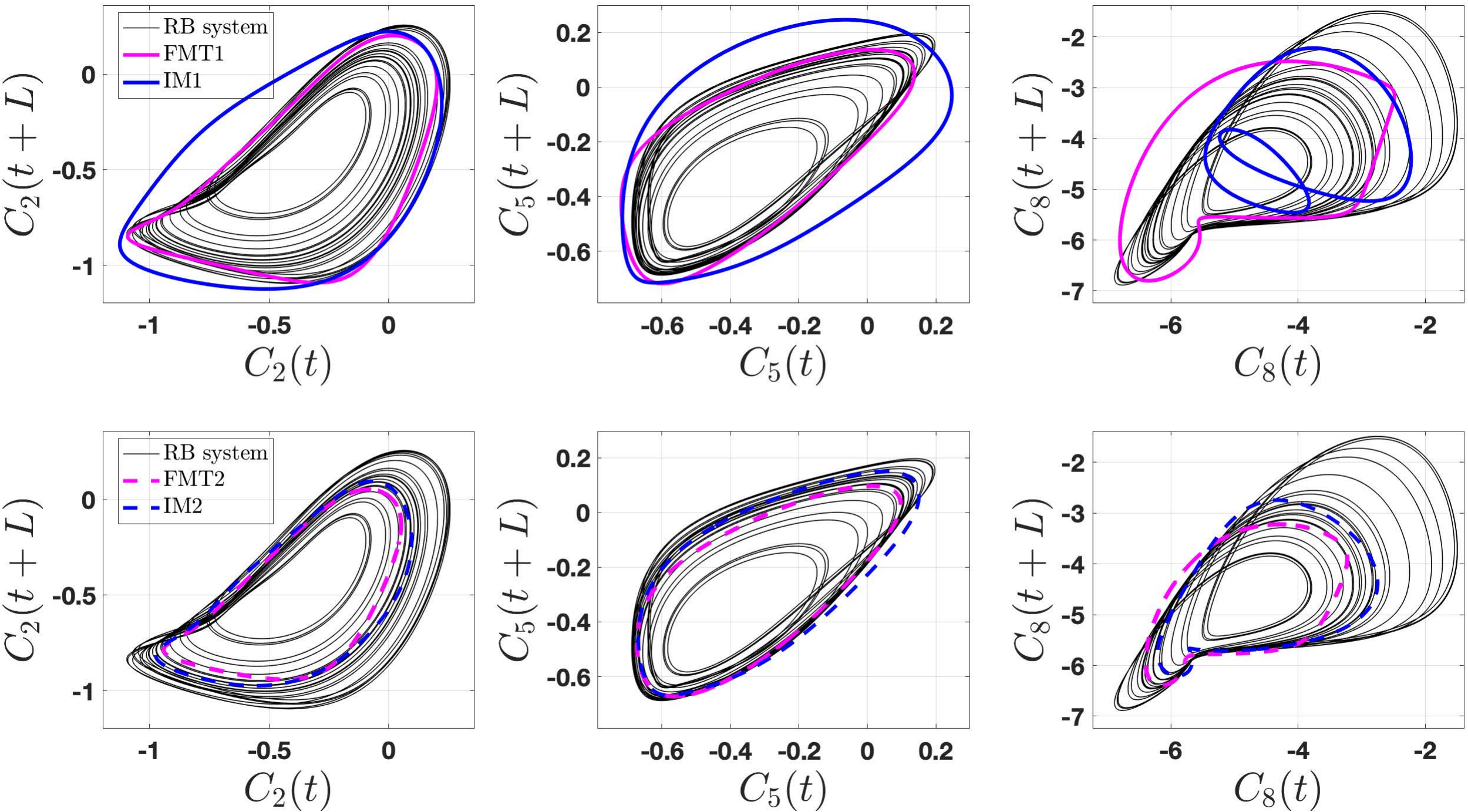}
\caption{The IM/FMT reduced systems are derived from  Eq.~\eqref{Eq_RBC_fluct} when the fluctuations are either taken with respect to the  mean state $\overline{\boldsymbol{C}}_{r}$ at $r=14.22$ (corresponding to the point $P_2$ in Fig.~\ref{Fig_RB_9D_combo}D) (bottom row) or with respect to the closest steady state to  $\overline{\boldsymbol{C}}_{r}$ (top row). Whatever the strategy retained, the IM/FMT reduced systems fail to predict the proper chaotic dynamics, predicting instead periodic dynamics for the dimension $m_c=5$ of the reduced state space than used for the OPM results shown in Fig.~\ref{Fig_RB_9D_combo}B}. \label{Fig_RB9D_h1_chaotic}
\end{figure}
%%%%%%%%%%%%%%%%%%%%%%%%%%%%%%%%%%%%%%%%%%%%%%%%%%%%%%
To the solution $z(t)=(z_1(t),\cdots,z_{m_c}(t))$ to this IM reduced system one then forms the approximation $\boldsymbol{C}^{\text{IM}}(t)$ of $\boldsymbol{C}(t)$ given by:
\be\label{Eq_CIM}
\boldsymbol{C}^{\text{IM}}(t) = \sum_{j=1}^{m_c} z_j(t) \boldsymbol{f}_j  + \sum_{n=m_c+1}^9 h_{n}( z(t)) \boldsymbol{f}_n +\overline{\boldsymbol{Y}}_{\! r}. 
\ee
%%%%%%%%%%%%%%%%%%%%%%
Similarly, we also test the performance of the reduced system when the local IM \eqref{Eq_h2} is replaced by the following FMT parameterization (see \eqref{Eq_FMTgen} in Remark \ref{Rem_FMT}):
\be \label{Eq_FMT}
h^{\mathrm{FMT}}_n = - \sum_{i, j = 1}^{m_c} \frac{\tilde{B}_{ij}^n}{\lambda_{n}} X_{i} X_{j}, \quad n = m_c+1, \cdots, 9.
\ee
Formulas such as \eqref{Eq_FMT} have been used in earlier studies relying on inertial manifolds for predicting higher-order bifurcations \cite{brown1991minimal}.

The predicted orbits obtained from reduced systems built either on the IM parameterization \eqref{Eq_h2} or the FMT one \eqref{Eq_FMT} are shown in the top row of Fig.~\ref{Fig_RB9D_h1_chaotic} as blue and magenta curves, respectively.   Both reduced systems lead to periodic dynamics and thus fail dramatically in predicting the chaotic dynamics. The small spectral gap issue mentioned in Remark \ref{Rem_smallgap} plays a role in explaining this failure but not only.
Another fundamental issue for the closure of chaotic dynamics lies in the difficulty to provide accurate parameterization of small-energy  but dynamically important  variables; see e.g.~\cite[Sec.~6]{CLM19_closure}. This issue is encountered here, as some of variables to parameterize for Experiment II contain only from $0.23\%$ to $2.1\%$ of the total energy.

Replacing $\overline{\boldsymbol{Y}}_{\! r}$ by the genuine mean state $\overline{\boldsymbol{C}}_{r}$ does not 
fix this issue as some of variables to parameterize contain still a small fraction of the total energy, from $0.36\%$ to $1.5\%$. The IM parameterization \eqref{Eq_h2} and the FMT one \eqref{Eq_FMT} whose coefficients are now determined from  the spectral elements of $\tilde{A}$ in which  $\overline{\boldsymbol{C}}_{r}$ replaces $\overline{\boldsymbol{Y}}_{\! r}$ in  \eqref{Eq_A_forFMT}, still fail in parameterizing accurately such small-energy variables; see bottom row of  Fig.~\ref{Fig_RB9D_h1_chaotic}\footnote{Note that here the FMT parameterization accounts for the forcing term (see \cite[Eq.~(4.40)]{CLM19_closure}) produced by  fluctuations calculated with respect to $\overline{\boldsymbol{C}}_{r}$}. 
By comparison, the OPM succeeds in parameterizing accurately these small-energy variables. 
For Experiment I, inaccurate predictions are also observed from the reduced systems built either on the IM parameterization \eqref{Eq_h2} or the FMT one \eqref{Eq_FMT} (not shown).

%%%%%%%%%%%%%%%%%%%%%%%%%%


\end{document}